\documentclass[11pt,onecolumn,a4paper]{elsart}
\usepackage{amstext,amssymb,amsmath,amscd}
\usepackage{graphics,graphicx,epsfig,subfigure,color,pstricks,rotating,dsfont,mathrsfs,theorem,algorithm,array,algorithmic,longtable,multirow,multicol}%

\newcommand{\parenth}[1]{\left(#1\right)}
\newcommand{\crochets}[1]{\left[#1\right]}
\newcommand{\abs}[1]{\left|#1\right|}

\newcommand{\EV}{{\mathscr{E}_{\mathbb{V}}}}
\newtheorem{theorem}{Theorem}

\newtheorem{proposition}{Proposition}
\newtheorem{lemme}{Lemma}

\def\QED{\mbox{\rule[0pt]{1.3ex}{1.3ex}}}
\newcommand{\proof}[1]{{\bf \itshape Proof #1: }}
\def\endproof{\hspace*{\fill}~\QED\par\endtrivlist\unskip}

\begin{document}

\begin{frontmatter}

\title{Simulation of Gegenbauer Processes using Wavelet Packets}

\author[qut]{J. J. Collet} and \author[greyc]{M. J. Fadili\corauthref{cor1}\thanksref{thank1}}

\corauth[cor1]{Dr. M.J.Fadili is an Associate Professor with the Image Processing Group GREYC CNRS UMR 6072 14050 Caen Cedex France. GREYC CNRS UMR 6072 - ENSICAEN 6, Bd du Mar\'echal Juin 14050 Caen Cedex France. Tel: +33-(0)31-45-29-20 Fax: +33-(0)31-45-26-98.}
\ead{Jalal.Fadili@greyc.ensicaen.fr}

\address[qut]{School of Economics and Finance, Queensland University of Technology, GP0 Box 2434, Brisbane QLD 4001, Australia.}
\address[greyc]{Image Processing Group GREYC CNRS UMR 6072 14050 Caen Cedex France.}

\thanks[thank1]{A part of this work was carried out while MJF was invited as a visiting scholar at the Queensland University of Technology (QUT). MJF would like to acknowledge support of a Faculty Research grant QUT, and thank JC for hospitality during the stay.}

\date{}

\maketitle

\begin{abstract}
In this paper, we study the synthesis of Gegenbauer processes using the wavelet packets transform. In order to
simulate a $1$-factor Gegenbauer process, we introduce an original algorithm, inspired by the one proposed by
Coifman and Wickerhauser \cite{CoiWic92}, to adaptively search for the best-ortho-basis in the wavelet packet
library where the covariance matrix of the transformed process is nearly diagonal. Our method clearly
outperforms the one recently proposed by \cite{Whi01}, is very fast, does not depend on the wavelet choice, and
is not very sensitive to the length of the time series. From these first results we propose an algorithm to
build bases to simulate $k$-factor Gegenbauer processes. Given its practical simplicity, we
feel the general practitioner will be attracted to our simulator. Finally we evaluate the approximation due to
the fact that we consider the wavelet packet coefficients as uncorrelated. An empirical study is carried out
which supports our results.
\end{abstract}

\begin{keyword}
Gegenbauer process, Wavelet packet transform, Best-basis, Autocovariance.
\end{keyword}

\end{frontmatter}


\section{Introduction}
The  simulation of long memory processes is an issue of a paramount importance in many statistical problems. In
the time domain, there exist different methods devoted to this task (see \cite{Ber94} for a non exhaustive
review of them). Alternative  efficient approaches, which operate in the frequency domain, were also proposed
(see \cite{Hosk84,DavHar87,Ber94}). More recently, owing to their scale-invariance property, wavelets have
since been widely adopted as a natural tool for analyzing and synthesizing $1/f$ long-memory processes. They
were demonstrated to provide almost Karhunen-Lo\`eve expansion of such processes \cite{Wor96}.

The simulation of fractional differenced Gaussian noise (fdGn) using discrete wavelet transform (DWT) has been
studied by \cite{MccWal96}. This kind of process is characterized by an unbounded power spectral density (PSD)
at zero. The proposed method relies on the fact that the DWT approximately decorrelates long memory processes (see e.g.
\cite{DerTew93,TewKim92,Wor96,Jen99,PerWal00}). The orthonormal wavelet decomposition "only" ensures
approximate decorrelation. The quality of this approximation has been widely assessed in \cite{Wor90,Flan92,Dijk94,TewKim92,Wor96,Jen99,Jen00} for a variety of $1/f$ long memory processes.

The DWT is only adapted to processes whose PSD is unbounded at the origin. Gegenbauer processes (sometimes also
called seasonal persistent processes) are also long memory processes and are characterized by an unbounded PSD.
The main difference with the fdGn processes is that the singularities of the PSD of the Gegenbauer processes can
be located at one or many frequencies in the Nyquist domain, not necessarily at the origin. Therefore, a natural
tool to analyze such processes appears to be the {\it wavelet packet transform} (WPT), which is a generalization of
the wavelet transform. The wavelets packets adaptively divide the frequency axis into separate dyadic intervals of
various sizes. They segment unconditionally, the frequency axis and are uniformly translated in time. Moreover, a
discrete time series of size $N$ is decomposed into more than $2^{N/2}$ wavelet packet (WP) bases. Among these bases,
one is a very good candidate to whiten the series and hence almost diagonalizes the covariance of the seasonal
process.

Recently, Mallat, Zhang and Papanicolaou \cite{MalZhaPap98}, and, following their work, Donoho, Mallat and von
Sachs \cite{DonMalVon98}, studied the idea of estimating the covariance of locally stationary processes by
approximating the covariance of the process by a covariance which is almost diagonal in a specially constructed
basis (cosine packets for their locally stationary processes) using an adaptation of Coifman-Wickerhauser (CW) best
ortho-basis algorithm. To some extent (given that we are interested in synthesis and they were in estimation issues), our work here can be seen as the spectral dual of theirs, since we are interested in studying the covariance of seasonal processes in the WP domain.

To the best of our knowledge, the simulation of the Gegenbauer process using the Discrete WPT (DWPT) has
been first studied in \cite{Whi01}. The DWPT creates a redundant collection of wavelet
coefficients at each level of the transform organized in a binary tree structure, equipped with a natural inheritance property. Different methods exist to determine the best candidate orthonormal basis. The author in \cite{Whi01} used a method which depends on both the location of the singularity and the wavelet used in the DWPT. To simulate realizations of a Gegenbauer process, once the basis is found, it then remains to apply the (inverse) DWPT using the same approximation as in \cite{MccWal96}.

This basis search method consists first in considering the square gain function of the wavelet filter associated with each
WP coefficient that is sufficiently small at the Gegenbauer frequency. Then a pruning of this family
is done to obtain the ortho-basis. The main advantage of this method is its simplicity. However
several points are still questionable and must be clarified. First the notion "sufficiently small" implies the
introduction of a threshold which seems to depend both on the wavelet used and the length of the simulated
series. No indication is given how to choose this threshold which remains awkward to control. Furthermore, it is not clear why
the basis should depend on the wavelet. Lastly, this method inherently leads to an
over-partitioning of the spectra which depends on the wavelet and the threshold considered (see e.g. Fig.\ref{basis_spec_freq} and Fig.4 in
\cite{Whi01}). Indeed, as the Gegenbauer process we consider here is stationary, it is known that the Karhunen-Lo\`eve basis is the Fourier basis. While over-partitioning, the approach of \cite{Whi01} inherently tries to approach the Fourier basis (more precisely it tends to select most of the atoms from the Shannon wavelet packets at the deepest level). Then, this makes wavelet packets machinery only of limited interest here. Furthermore, many important statistical tasks involving Gegenbauer processes would seriously suffer from such an over-partitioning, e.g. maximum likelihood estimation, resampling-based inference, to cite only a few examples. 

To alleviate these intricacies, our belief is that it should more beneficial to build, for each Gegenbauer process, an unique valid (almost whitening) basis for all wavelets with a reduced number of packets. This basis should only depend on the Gegenbauer frequencies, but not on the long memory parameters nor on the wavelet used. The rationale behind these claims can be supported by different arguments. Indeed, wavelets are now widespread as almost-diagonalizing expansion for $1/f$ processes, no matter what the long memory parameter and the wavelet are. Although, the latter parameters clearly influence the quality of the decorrelation as was widely proven \cite{Wor96}. Our goal is then to mimic this behavior by extending and generalizing the aforementioned properties to Gegenbauer processes within the WP framework, with the desirable properties that (i) our basis tends to the dyadic wavelet basis (the 1-band WP) when the singularity frequency tends to $0$, and (ii) the provided basis should have a limited number of packets. To get a gist of the latter property, we can say that we are seeking a basis (and the corresponding tree) which attains the minimal diagonalization error penalized by the complexity of the tree in terms of the number of packets (i.e. number of leaves of the tree) involved in the dyadic partition of the spectral axis provided by the selected WP basis. See \cite[Sec. 13]{DonMalVon98} and \cite{Donoho97} for a more detailed discussion of complexity penalized estimation and its relation to best-ortho-basis.
 
In this paper, we propose an alternative efficient way to determine the appropriate basis for the simulation of
$1$-factor Gegenbauer process, that we extend to the simulation of $k$-factor Gegenbauer process. To find this basis, We propose an algorithm which is an adaptation of the best-basis search algorithm of \cite{CoiWic92}. The main property of this algorithm is that it provides us with a (unique) basis that only depends on the location of the Gegenbauer frequency, unlike the construction method of \cite{Whi01} which provides bases depending both on the location of the singularity and the wavelet. To point out the role played by the wavelet used, we will study the decorrelation properties of the WP coefficients of a Gegenbauer process when it is
expressed in this basis. In particular, the influence of the wavelet regularity, the long memory parameter and the location of the singularity on the decorrelation decay speed will be established.

The organization of this paper is as follows. After some preliminaries and notations related to the WPT theory (Section \ref{sec:21}) and to the Gegenbauer process (Section \ref{sec:22}) are introduced, we will define the best-basis search algorithm and the cost function we propose (Sections \ref{sec:31}-\ref{sec:32}). Theoretical support to this cost function is also supplied. We then develop an algorithm to build an appropriate basis to simulate $1$-factor Gegenbauer process (Section \ref{sec:33}). This method will then be extended to $k$-factor processes (Section \ref{sec:34}). Theoretical evaluation of the approximation quality due to the fact that we consider the WP coefficient as uncorrelated is studied in Section \ref{sec:4}. A simulation study is finally conducted to illustrate and discuss our results (Section \ref{sec:5}).

\section{Preliminaries}
\subsection{The wavelet packet transform}
\label{sec:21}
Wavelet packets were introduced by Coifman, Meyer and Wickerhauser \cite{CoiMeyWic92}, by generalizing
the link between multi-resolution approximations and wavelets. Let the sequence of functions defined recursively
as follows:
\begin{align}
\psi_{j+1}^{2p}(t) & =\sum_{n=-\infty}^{\infty}h(n)\psi_j^p(t-2^jn) \\
\psi_{j+1}^{2p+1}(t) & =\sum_{n=-\infty}^{\infty}g(n)\psi_j^p(t-2^jn)
\end{align}
for $j\in\mathbb{N}$ and $p=0,\dots,2^j-1$, where $h$ and $g$ are the conjugate pair of quadrature mirror
filters (QMF). At the first scale, the functions $\psi_0$ and $\psi_1$ can be respectively identified with the
father and the mother wavelets $\phi$ and $\psi_1$ with the classical properties (among others):
\begin{equation}
\int \phi(t) = 1, \int \psi(t) = 0
\end{equation}

The collection of translated, dilated and normalized functions
$\psi^{p,n}_{j}\overset{\mathrm{def}}{=}2^{-j/2}\psi_p(2^{-j}t-n)$ makes up what we call the (multi-scale) wavelet
packets associated to the QMFs $h$ and $g$. $j\in\mathbb{N}$ is the scale index, $p=0,\ldots,2^j-1$ can be
identified with a frequency index and $k$ is the position index. It has been proved (see e.g. \cite{Wic94}) that if $\{\psi_j^{p,n}\}_{n\in\mathbb{Z}}$ is an orthonormal basis of a space $\mathbf{V}_j$, then the
family $\{\psi_{j+1}^{2p,n},\psi_{j+1}^{2p+1,n}\}_{n\in\mathbb{Z}}$ is also an orthonormal basis of $\mathbf{V}_j$.

The recursive splitting of vector spaces is represented in a binary tree. To each node $(j,p)$, with
$j\in\mathbb{N}$ and $p=0,\dots,2^j-1$, we associate a space $\mathbf{V}_j^p$ with the orthonormal basis
$\{\psi_{j}^{p}(t-2^{j}n)\}_{n\in\mathbb{Z}}$. As the splitting relations creates two orthogonal basis, it is
obvious that $\mathbf{V}_j^p=\mathbf{V}_{j+1}^{2p}\oplus\mathbf{V}_{j+1}^{2p+1}$.

The WP representation is overcomplete. That is, there are many subsets of wavelet packets which
constitute orthonormal bases for the original space $\mathbf{V_0}$ (typically more than $2^{2^{J-1}}$ for a
binary tree of depth $J$). While they form a large library, these bases can be easily organized in a binary tree
and efficiently searched for extreme points of certain cost functions, see \cite{CoiWic92} for
details. Such a search algorithm and associated cost function are at the heart of this paper.

In the following we call the collection $\mathcal{B}=\{\psi_j^{p,n}\}_{(j,p) \in \mathcal{T}, n \in \mathbb{Z}}$ the basis of
$\mathbf{L^2}(\mathbb{R})$, and the tree $\mathcal{T}$ for which the collection
of nodes $(j,p)$ are the leaves, the associated tree.


Given a basis $\mathcal{B}$ and its associated tree $\mathcal{T}$ it is possible to decompose any function $x$
of $\mathbf{L^2}(\mathbb{R})$ in $\mathcal{B}$. At each node $(j,p) \in \mathcal{T}$,
the WP coefficients $W_{j}^{p}(n)$ of $x$ in the subspace $\mathbf{V}_{j}^{p}$ at position
$n$ are given by the inner product:
\begin{equation}
W_{j}^{p}(n)=\int\psi_{j}^{p}(t-2^{j}n)x(t)dt.
\end{equation}
For a discrete signal of $N$ equally-spaced samples, the DWPT is calculated using a fast
filter bank algorithm that requires $O(N\log N)$ operations. The interested reader may refer to the books of
Mallat \cite{Mal98} and Wickerhauser \cite{Wic94} for more details about the DWPT.

\subsection{Gegenbauer process}
\label{sec:22}
The $k$-factor Gegenbauer process is a $1/f$-type process introduced in \cite{GraZhaWoo89,Gray98}.
The PSD $f$ of a such process $(X_t)_t$ is given by for all $|\lambda|\leq1/2$
\begin{equation}
\label{GG_spect_dens}
f(\lambda)=\frac{\sigma_\varepsilon^2}{2\pi}\prod_{i=1}^{k}\parenth{4\parenth{\cos2\pi\lambda-\cos2\pi\nu_i}^2}^{-d_i}
\end{equation}
where $k$ is a finite integer and $0<d_i<1/2$ if $0<|\nu_i|<1/2$ and $0<d_i<1/4$ if $|\nu_i|=0$ for
$i=1,\dots,k$. The parameter $d_i$ and $\nu_i$ are respectively called the memory parameter and the Gegenbauer
frequency. The $k$-factor Gegenbauer process is a generalization of the fractionally differenced Gaussian white
noise process (see \cite{Hos81} and \cite{GraJoy80}) in the sense that
the PSD is unbounded at $k$ different frequencies not necessary located in $0$.

The Gegenbauer process $(X_t)_t$ is related to a white noise process $(\varepsilon_t)_t$ with mean $0$ and
variance $\sigma_\varepsilon^2$ through the relationship:
\begin{equation}
\label{def2} \prod_{i=1}^{k} (I-2\nu_iB+B^2)^{d_i}X_t=\varepsilon_t,
\end{equation}
\noindent where $BX_t=X_{t-1}$ and $\eta_i=\cos2\pi\nu_i$.\\
The main characteristic of the Gegenbauer processes in the time domain is the slow decay of autocovariance
function. In the case of a $1$-factor Gegenbauer process, Gray {\it et al.} \cite{GraZhaWoo89} and then
Chung \cite{Chu96} proved the asymptotic behavior of the autocovariance function:
\begin{equation}
\rho(h)\sim h^{2d-1}\cos(2\pi\nu h) \quad \textrm{as} \quad h\rightarrow\infty.
\end{equation}
The next section is devoted to the construction of the best basis diagonalizing the covariance of a $N$-sample
realization of a Gegenbauer process with the convention $N=2^J$.

\section{Simulation of Gegenbauer processes}
This section is composed of two parts. The first one is devoted to the simulation procedure in the general case:
no assumption is made concerning the basis, except that we have an appropriate basis. The second part concerns
the construction of this appropriate basis.

\subsection{Simulation procedure}
\label{sec:31}
Here we present the procedure to simulate a Gegenbauer process. Assume we would like to simulate a $k$-factor
Gegenbauer process, $(X_t)_t$, with PSD $f$ as defined in (\ref{GG_spect_dens}), with Gegenbauer
frequencies $(\lambda_1,\dots,\lambda_k)$ and long memory parameters $(d_1,\dots,d_k)$. The length of the
realization will be $N=2^J$.

We define the band-pass variance $\beta^2_{j,p}$ in the frequency interval
$I^p_j=[\frac{p}{2^{j+1}},\frac{p+1}{2^{j+1}}]$ by:
\begin{equation}
\label{Eq:Bjp}
\beta^2_{j,p}=2\int_{\frac{p}{2^{j+1}}}^{\frac{p+1}{2^{j+1}}}f(\lambda)d\lambda
\end{equation}
As in \cite{MccWal96} and \cite{Whi01}, we assume that the PSD in each frequency interval $I^p_j$, for which the couple $(j,p)$ is a leaf of the tree $\mathcal{T}$ associated to the basis $\mathcal{B}$, is constant and equal to $\sigma^2_{j,p}$. Then, the band-pass variance is (approximately) equal to:
\begin{equation}
\beta^2_{j,p}=2\int_{\frac{p}{2^{j+1}}}^{\frac{p+1}{2^{j+1}}}\sigma^2_{j,p}d\lambda=2^{-j}\sigma^2_{j,p}
\end{equation}
Thus the variance of each WP coefficient is given by $\mathbb{V}[W_j^p(n)]=\sigma^2_{j,p}=2^j\beta^2_{j,p}$, $\mathbb{V}$ is the variance operator. To simulate $N$ observations of a Gegenbauer process $(X_t)_{t=1,\ldots,N}$ with PSD $f$, we use the following
procedure:

\begin{algorithmic}[1]
\STATE Given an appropriate basis $\mathcal{B}$ and its associated tree $\mathcal{T}$, calculate the band-pass variances $\beta^2_{j,p}$, $(j,p)\in\mathcal{T}$ as in (\ref{Eq:Bjp});
\STATE For each $(j,p)\in\mathcal{T}$, generate $2^{J-j}$ realizations of 
$W_j^p(n)$, an independent Gaussian random variable with zero mean and variance equal to
$\sigma^2_{j,p}$; 
\STATE Organize the WP coefficients $W_j^p(n)$, for $(j,p)\in\mathcal{T}$ and $n=1,\dots,2^{J-j}$, in a
vector $\mathbf{W}_\mathcal{B}$, and apply the the inverse DWPT to obtain the
observation vector $\mathbf{X}=(X_1,\dots,X_N)^T$.
\end{algorithmic}

In the following subsection, we examine the construction of what we term an appropriate basis $\mathcal{B}$.

\subsection{Best-basis construction algorithm}
\label{sec:32}
\subsubsection{Approximate Diagonalization in a Best-Ortho-basis} 
\label{sec:BoBDon} 
Let $(X_t)_t$ be a stationary Gegenbauer process and $\Gamma$ its covariance matrix. Let $\gamma_{i,j}\crochets{{\mathcal B}}$ the
entries of $\Gamma\crochets{\mathcal B}$; the covariance matrix of the
coordinates $\mathbf{W}_\mathcal{B}$ of $(X_t)_t$ in the ortho-basis ${\mathcal B}$. One can define diagonalization as an optimization of
the functional \cite{DonMalVon98}:
\begin{equation}
\max_{{\mathcal B}} {\mathcal E}({\mathcal B}) = \max_{{\mathcal B}} \sum_i e(\gamma_{ii}[{\mathcal B}])
\end{equation}
where $e$ is taken as a strictly convex cost function. In practice, the optimization formulation of diagonalization is not
widely used, presumably because it generally does not help in computing diagonalizations. Optimization of an
arbitrary objective ${\mathcal E}$ over finite libraries of orthogonal bases - the cosine packets library and
the wavelet packets library - is not a problem with good algorithmic solutions. Wickerhauser \cite{Wic91}
suggested applying these libraries in problems related to covariance estimation. He proposed the notion of
selecting a "best basis" for representing a covariance by optimization of the "entropy functional" $e_H(\gamma) = -\log \gamma$ over all bases in a restricted library. Authors in \cite{MalZhaPap98}, developed a proposal
which uses the specific choice $e_2(\gamma) = \gamma^2$.


In the Wickerhauser formulation, one is optimizing over a finite library and there will not generally be a basis
in this library which exactly diagonalizes $\Gamma$. Then different strictly convex functions $e(\gamma)$ may end up
picking different bases. For example, the quadratic cost function $e_2$ has a special interpretation in this context as it
leads to a basis which best diagonalizes $\Gamma$ in a least-squares sense \cite{DonMalVon98}, and is closely related to the Hilbert-Schmidt (HS) norm of the diagonalization error. Similarly, the -log "entropy functional" is connected to the Kullback-Leibler divergence \cite{DonMalVon98}. Even if the approach developed in \cite{MalZhaPap98,DonMalVon98} was
specialized to the case of $e_2$, it is not really tied to the specific entropy measure; other additive convex
measures can be accommodated such as the $l_\alpha$ norm $\alpha>2$ or the neg-entropy, and the CW proposal makes equally sense. This was the starting point of our work.


\subsubsection{Proposed Algorithm} 
The optimization problem of ${\mathcal E}$ over bases can be re-expressed as an optimization over trees, as follows. Set $\EV[W_{j}^{p}]=\sum_{n_j} e\parenth{\mathbb{V}\crochets{W_{j}^{p}(n_j)}}$. Then as $\sum_{\mu \in \mathcal{B}}=\sum_{(j,p) \in \mathcal{T}}\sum_{n_j}$, one is actually trying to optimize:
\begin{equation}
\sum_{(j,p) \in \mathcal{T}} \EV[W_{j}^{p}]
\end{equation}
over all recursive dyadic partitions of the spectral axis. The best basis $\mathcal{B}$ is then the one that maximizes some measure of the wavelet packets variances, among all the bases that can be constructed from the tree-structured library. The construction of the best basis can be accomplished efficiently using the recursive bottom-up CW algorithm defined by \cite{CoiWic92}:
\begin{equation}
\label{algo_CW} \mathcal{B}_j^p=\begin{cases}
            \mathcal{B}_{j+1}^{2p}\cup\mathcal{B}_{j+1}^{2p+1} & \textrm{if} \quad \EV[W_{j+1}^{2p}] +
            \EV[W_{j+1}^{2p+1}]>\EV[W_{j}^{p}],\\
            \mathcal{B}_{j}^{p} & \textrm{if} \quad \EV[W_{j+1}^{2p}] +
            \EV[W_{j+1}^{2p+1}]\leq \EV[W_{j}^{p}].
        \end{cases}
\end{equation}

The chosen criterion lies on the comparison between some measure of WP coefficients variances at the children nodes and their parents. Beside the fact that these variances, and then the basis, will depend on the long-memory parameter (or even the wavelet), there is another even more important reason that prevents from a crude use of such a search algorithm with the cost functions that we defined above (such as $e_2$). Indeed, the band-pass variance of any node is equal to the sum of those of its children. Hence, it is not a difficult matter to check that any strictly convex cost functional such as those specified above, e.g. $e_2$ or $e_H$, will systematically provide the basis corresponding to the finest partition of the spectral axis, which is clearly the worst in terms of complexity (i.e. number of wavelet packets). Again, this makes the wavelet packets machinery only of limited interest.


Therefore, motivated by the above discussion, we were led to define, for the wavelet packets $W_{j+1}^{2p}$ and $W_{j+1}^{2p+1}$, a new type of WP variance cost measure as follows:
\begin{equation}
  \EV[W_{j+1}^{2p}]=\begin{cases}
                    0 & \textrm{if} \quad \mathbb{V}[W_{j+1}^{2p}]\leq A_0\mathbb{V}[W_{j+1}^{2p+1}]\\
                    \mathbb{V}[W_{j+1}^{2p}] & \textrm{otherwise.}
                      \end{cases}
\end{equation}
\begin{equation}
  \EV[W_{j+1}^{2p+1}]=\begin{cases}
                    0 & \textrm{if} \quad \mathbb{V}[W_{j+1}^{2p+1}]\leq A_0\mathbb{V}[W_{j+1}^{2p}]\\
                    \mathbb{V}[W_{j+1}^{2p+1}] & \textrm{otherwise.}
                      \end{cases}
\end{equation}
where $A_0$ is a fixed positive constant (its value will depend for instance on the singularity frequency and will
be given in the proof of Proposition \ref{Prop_1_fact}).

In the following, when we write (with a slight abuse of notation) that $\EV[W_{j+1}^{2p}]=0$ or
$\EV[W_{j+1}^{2p+1}]=0$, it will mean respectively that there exists a constant $A_0<1$ such that $\mathbb{V}[W_{j+1}^{2p}]\leq
A_0\mathbb{V}[W_{j+1}^{2p+1}]$ or $\mathbb{V}[W_{j+1}^{2p+1}]\leq A_0\mathbb{V}[W_{j+1}^{2p}]$. In these cases
we will also use respectively the notations,
$$\mathbb{V}[W_{j+1}^{2p}]\ll\mathbb{V}[W_{j+1}^{2p+1}]\ \ \ \ \ \textrm{and}\ \ \ \ \ \
\mathbb{V}[W_{j+1}^{2p+1}]\ll\mathbb{V}[W_{j+1}^{2p}]$$

Using the criterion defined above, algorithm (\ref{algo_CW}) becomes\footnote{Strictly speaking, this is no longer a CW algorithm.}: 
\begin{equation} 
\label{new_algo_CW}
\mathcal{B}_j^p=\begin{cases}
         \mathcal{B}_{j+1}^{2p}\cup\mathcal{B}_{j+1}^{2p+1}, & \textrm{if} \quad
         \EV[W_{j+1}^{2p}]=0 ~\textrm{or}~\EV[W_{j+1}^{2p+1}]=0,\\
         \mathcal{B}_{j}^{p}, & \textrm{otherwise.}
        \end{cases}
\end{equation}
In the following, we use this algorithm to build the best-ortho-basis for a Gegenbauer process.

\subsection{The $1$-factor case} 
\label{sec:33}
It is natural to build the best basis according to the shape of the PSD of our process. More precisely, the basis is a function of the location of the singularities. It means that in the case of $1$-factor Gegenbauer process, the basis depends directly on the value of the Gegenbauer frequency. Using the notations defined in the previous section, the recursive construction is summarized in the
following proposition.
\begin{proposition}
  \label{Prop_1_fact}
  If $(X_t)_t$ is a stationary $1$-factor Gegenbauer process, with parameters $(d,\nu,\sigma)$
  then, at node $(j,p)$, if the frequency $\nu$ is in the interval $I_{j}^{p}=[\frac{p}{2^j},\frac{p+1}{2^j}[$,
  then:
  \[
  \EV[W_{j+1}^{2p}]=0\ \ \ \ \ \textrm{or}\ \ \ \ \ \EV[W_{j+1}^{2p+1}]=0,
  \]
  and consequently for algorithm (\ref{new_algo_CW}):
  \[
  \mathcal{B}_j^p=\mathcal{B}_{j+1}^{2p}\cup\mathcal{B}_{j+1}^{2p+1}.
  \]
  Furthermore, if the frequency $\nu$ is in the closure of the intervals $I_{j+1}^{2p}$ and $I_{j+1}^{2p+1}$,
  then:
  \[
  \EV[W_{j+2}^{4p+1}]=0\ \ \ \ \ \textrm{and}\ \ \ \ \ \EV[W_{j+2}^{4p+2}]=0,
  \]
  and consequently for algorithm (\ref{new_algo_CW}):
  \[
  \mathcal{B}_j^p=\mathcal{B}_{j+2}^{4p}\cup\mathcal{B}_{j+2}^{4p+1}\cup\mathcal{B}_{j+2}^{4p+2}\cup
  \mathcal{B}_{j+2}^{4p+3}.
  \]
\end{proposition}

Proof:\rm\ {\it See Appendix A}.

To construct the best-ortho-basis of a $1$-factor Gegenbauer process, we propose Algorithm \ref{algo:1} which proceeds according to
Proposition \ref{Prop_1_fact} and the aggregation relation defined in (\ref{new_algo_CW}).
\begin{algorithm}
\caption{$1$-factor Best-Basis Search Algorithm} 
\label{algo:1}
\begin{algorithmic}[1]
\label{1factoralgo} 
\REQUIRE{A Gegenbauer frequency $\nu$ and sample size $N=2^J$,}

\vspace{0.25cm} 
{\it \underline{\bf Initialization}} 
\FOR{$j=0,\dots,J$ and $p=0,\dots,2^j-1$}
  \STATE    $Tree(j,p)=0$.
\ENDFOR

\vspace{0.25cm} 
{\it \underline{\bf Main Loop}} 
\FOR{$j=1,\dots,J$}
  \FOR{$p=0,2,\dots,2^j-2$}
    \IF{$\nu\in[p/2^{j+1},(p+1)/2^{j+1}]$}
      \STATE $Tree(j,p+1)=1$
    \ENDIF
    \IF{$\nu\in[(p+1)/2^{j+1},(p+2)/2^{j+1}]$}
      \STATE $Tree(j,p)=1$
    \ENDIF
  \ENDFOR
\ENDFOR

\vspace{0.25cm} 
{\it \underline{\bf Pruning}} 
\FOR{$j=1,\dots,J$}
  \FOR{$p=0,2,\dots,2^j-2$}
    \IF{$\nu\in[p/2^{j+1},(p+1)/2^{j+1}]$ and $\mathop{\max}_{r=1,\dots,J-j-1; s=0,\dots,2^r-1} Tree(j+r,2^{r}p+s)>0$}
      \STATE $Tree(j,p)=0$
    \ENDIF
  \ENDFOR
\ENDFOR
\end{algorithmic}
\end{algorithm}
This algorithm is decomposed into two mains loops. The first one builds a family where the best-ortho-basis is
included. The second loop is a pruning of the family to obtain the best-ortho-basis. This second loop
corresponds to the second part of Proposition \ref{Prop_1_fact}.

\begin{figure}[H]
  \begin{center}
    \includegraphics[width=0.8\textwidth,height=6.9cm]{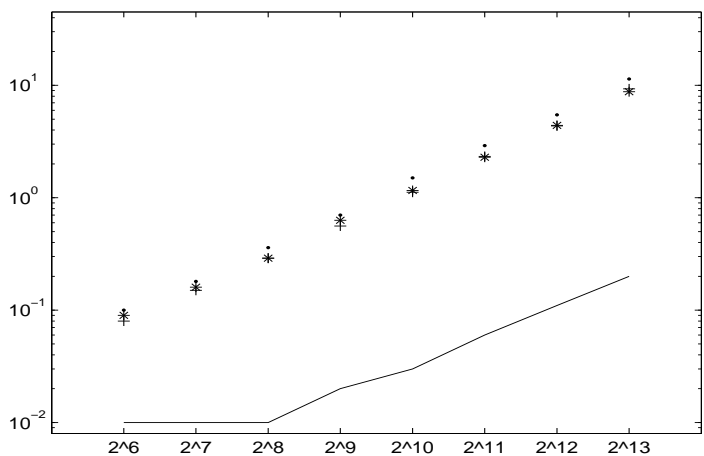}
    \caption{Computation time (in seconds) to build the best-ortho-basis. Solid line (our approach). The symbols $'+'$, $'\times'$ and $'.'$ correspond to the computation time using the method of \cite{Whi01} in the case of respectively $'db10'$ (Daubechies wavelet $q=10$), $'sym10'$ (Symmlet $q=10$) and $'coif5'$ (Coiflet $q=10$).}
    \label{Tps_construct_basis}
  \end{center}
\end{figure}

The algorithm we propose is very fast involving only simple comparisons, and it does not require the calculation of variances of WP
coefficients. To illustrate the computational speed of our algorithm we provide in Fig.\ref{Tps_construct_basis} some
computation times to build bases using our method and the method of \cite{Whi01}\footnote{The experiments were run under the R environment on a 2.4GHz PC with 512MB RAM}. In this example, we are only interested in the time needed to build the basis. These bases are built to simulate Gegenbauer process with a singularity located at $1/12$ and length equal to $2^J$, with $J=6,\dots,13$. The solid line corresponds to the computation time of the algorithm we propose. The symbols $'+'$, $'\times'$ and $'.'$ correspond to the computation time using the method of \cite{Whi01} in the case of respectively $'db10'$ (Daubechies wavelet with $q=10$ vanishing moments), $'sym10'$ (Symmlet $q=10$) and $'coif5'$ (Coiflet $q=10$). In every case the computation time increases with the length of the process. However this time increases always much faster for \cite{Whi01} than for the current method (the ratio of computation times is 10 to 300 times larger for the competitor method for series of length $64$ to $8192$). Typically, for a 8192-sample series, it takes 100 ms to our algorithm to find the best basis while
\cite{Whi01} algorithm requires 30 s.

\subsubsection*{Examples}
\begin{figure}[H]
  \centering
  \hspace*{-1cm}
  \includegraphics[width=1.15\textwidth]{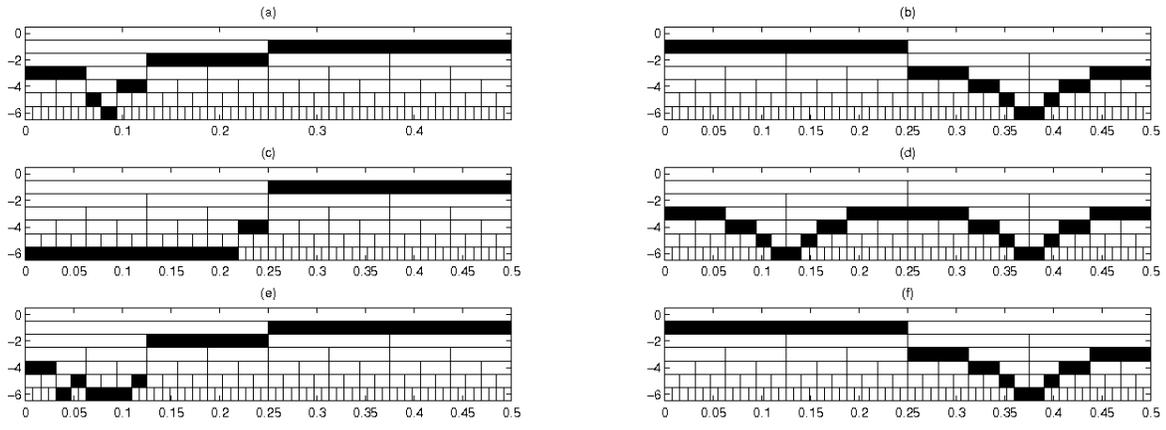}
   \caption{Best-orthos-basis for a Gegenbauer process, with $\nu=1/12$ (first column) and $\nu=0.375$ (second column). (a) Our basis $\nu=1/12$, (b) our basis $\nu=1/5$. Basis of \cite{Whi01} for $\nu=1/12$ with (c) $'db3'$ and (e) $'coif5'$ wavelets, and similarly for $\nu=0.375$ (d)-(f). Black rectangles correspond to the leaves of the binary tree, and then to the partition of the spectral axis.}
   \label{basis_spec_freq}
\end{figure}

We give two examples of construction of bases. Fig.\ref{basis_spec_freq}.(a) depicts the basis built using the first part of Algorithm \ref{algo:1}, to simulate a stationary Gegenbauer process with frequency $\nu=1/12$. Fig.\ref{basis_spec_freq}.(b) shows the basis constructed in the case of a stationary Gegenbauer process with $\nu=0.375$. The last case corresponds to the second situation of Proposition \ref{Prop_1_fact}. One may remark that unlike the first case where the tree has at least one leaf at each scale, in this second case, because of the particular value of the Gegenbauer frequency, there exists scale for which the tree has no leaf (see scale $j=2$). For comparative purposes, observe that the basis provided by the approach of \cite{Whi01} (with a threshold 0.01) is highly dependent of the wavelet choice. For example in Fig.\ref{basis_spec_freq}.(c) ($'db3'$), one cannot have an idea of the singularity location. In Fig.\ref{basis_spec_freq}.(d) ($'coif5'$), two singularties are apparent while only one is relevant. In both last cases, the basis is clearly over-partitioned.

\subsection{The $k$-factor case} 
\label{sec:34}
In this section we are interested in the general case: the construction of the
appropriate basis to simulate a $k$-factor Gegenbauer process. To achieve this goal, let us consider $(X_t^1)_t$
and $(X_t^2)_t$ as respectively a $(k-1)$-factor and a $1$-factor Gegenbauer processes. We denote
$(d_1,\nu_1,\dots,d_{k-1}, \nu_{k-1})$ and $(d_{k},\nu_{k})$ the parameters of $(X_t^1)_t$ and $(X_t^2)_t$. Let
$\mathcal{B}_1$ and $\mathcal{B}_2$ the best-ortho-bases of $(X_t^1)_t$ and $(X_t^2)_t$. We denote respectively
$\mathcal{T}_1$ and $\mathcal{T}_2$ the
trees associated with the bases $\mathcal{B}_1$ and $\mathcal{B}_2$.

Let $(X_t)_t$ be a $k$-factor Gegenbauer process with parameters $(d_1,\nu_1,\dots,d_{k},\nu_{k})$. We denote
$\mathcal{B}$ the appropriate basis and $\mathcal{T}$ the associated tree. Let $\mathcal{B}'$ be the family
equal to the union of the bases $\mathcal{B}_1$ and $\mathcal{B}_2$ and let $\mathcal{T}'$ be the associated
tree. We are now ready to state the following,

\begin{proposition}
  \label{Prop_k_fact}
  Under the previous assumptions:
  \begin{enumerate}
    \item $\mathcal{B}\subset\mathcal{B}'$
    \item Let $(j,p)$ be node in the tree $\mathcal{T}'$ such that there exists $r^*=1,\dots,J-j$ and
    $s^*=0,\dots,2^{r^*}-1$ such that $(j+r^*,2^{r^*}p+s)$ is also in the tree $\mathcal{T}'$. Then:
    \[
    (j,p)\not\in\mathcal{T}\ \ \ \ \ \textrm{and}\ \ \ \ \ (j+r^*,2^{r^*}p+s)\in\mathcal{T}
    \]
  \end{enumerate}
\end{proposition}

Proof:\rm\ {\it See Appendix A}.

According to this last proposition, the best-ortho-basis of a $k$-factor Gegenbauer process may be built using
$k$ well chosen best-ortho-bases of $1$-factor Gegenbauer processes. 
The steps outlined in Algorithm \ref{algo:2} allow to build the appropriate basis to simulate a $k$-factor Gegenbauer process. This
procedure lies on Algorithm \ref{1factoralgo} and results given in Proposition \ref{Prop_k_fact}.

\begin{algorithm}
\caption{$k$-factor Best-Basis Search Algorithm}
\label{algo:2}
\begin{algorithmic}[1]
\label{kfactoralgo} \REQUIRE{Gegenbauer frequencies $\nu_i$ and sample size $N=2^J$,}

\vspace{0.25cm} 
{\it \underline{\bf Initialization}} 
\FOR{Each Gegenbauer frequency $\nu_i, i=1,\dots,k$}
   \STATE Construct the best-ortho-basis $\mathcal{B}_i$ and associated tree $Tree_i$ using Algorithm \ref{algo:1}.
\ENDFOR 
\STATE $Tree=\cup_{i=1}^{k}Tree_i$ (implemented using e.g. the logical OR operator under R or Matlab).

\vspace{0.25cm} 
{\it \underline{\bf Pruning}} 
\FOR{$j=1,\dots,J$}
  \FOR{$p=0,2,\dots,2^j-2$}
    \IF{$Tree(j,p)=1$ and $\mathop{\max}_{r=1,\dots,J-j-1; s=0,\dots,2^r-1} Tree(j+r,2^{r}p+s)>0$}
      \STATE $Tree(j,p)=0$
    \ENDIF
  \ENDFOR
\ENDFOR
\end{algorithmic}
\end{algorithm}

\subsubsection*{Example}
\begin{figure}[H]
\vspace{-0.65cm}
  \begin{center}
      \includegraphics[width=\textwidth,height=9cm]{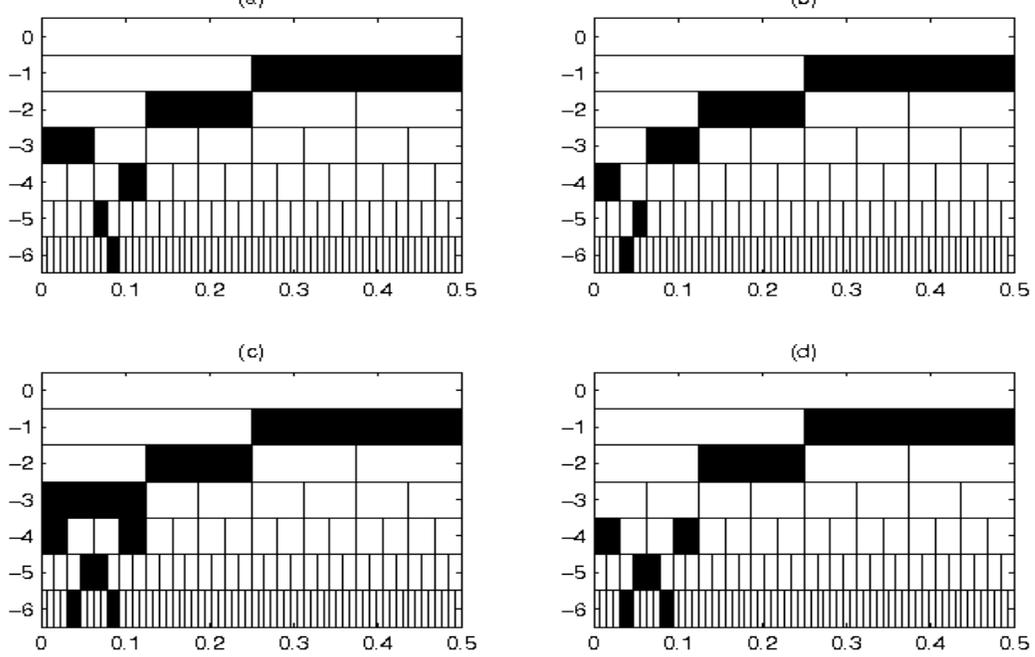}
  \end{center}
  \caption{(a) Best-ortho-basis $\mathcal{B}_1$. (b) Best-ortho-basis $\mathcal{B}_2$. (c) Union of bases $\mathcal{B}_1$ and $\mathcal{B}_2$. (d) Best-ortho-basis for the two factor Gegenbauer process $(X_t)_t$.}
  \label{basis_BoB}
\end{figure}

Here, we give an example of construction of the best-ortho-basis for a $2$-factor
Gegenbauer process $(X_t)_t$ with Gegenbauer frequencies $1/12$ and $1/24$. Fig.\ref{basis_BoB}.(a) and
\ref{basis_BoB}.(b) show the best-ortho-bases $\mathcal{B}_1$ and $\mathcal{B}_2$ of the
processes $(X^1_t)_t$ and $(X^2_t)_t$ (see the previous section for construction of these bases). The family $\mathcal{B}^*$ equal to $\mathcal{B}_1\cup\mathcal{B}_2$ is given in Fig.\ref{basis_BoB}.(c). This family is not a basis, the intersections between its elements are not always empty, e.g. at depth $j=3$, the elements at $p=0$ and $p=1$ should not be considered as elements of the best-ortho-basis and must be pruned away. This is accomplished using the methodology developed above, and an appropriate basis for the process $(X_t)_t$ is obtained as represented in Fig.\ref{basis_BoB}.(d).

\subsection{Back to the original CW algorithm}
Our aim here is to shed light on our best-basis search algorithm by relating it to the original CW one. More precisely, we shall give an additive cost functional, which can be used within the CW algorithm, that is closely linked to our proposal in (\ref{new_algo_CW}). Basically, the aggregation relation (\ref{new_algo_CW}) can be thought of as a rule which at each level, enforces dyadic splitting of an interval $I^p_j$ if (and only if) the WP variance inside that interval is above a threshold, and the node corresponding to this interval is marked as a branch (non-terminal). Otherwise the interval is kept intact (marked as a leaf) and the children nodes in the tree are pruned away. Doing so, this procedure implicitly tries to track the packets that contain the singularities of the process. Hence, motivated by these observations, an additive variance cost functional satisfying this aggregation rule can be defined as:
\begin{equation}
\label{Eq:costthresh1}
\EV[W^p_j]=\beta^2_{j,p} \mathds{1}\parenth{\beta^2_{j,p} \geq \delta}
\end{equation}
where $\beta^2_{j,p}$ is the band-pass variance as before and $\delta$ is a strictly positive threshold. The original recursive (bottom-up) CW could then be used to minimize such a cost functional (termed as "Number above a threshold" functional in Wickerhauser book \cite{Wic94}). Unfortunately, the threshold remains an important issue to fix, and depends jointly on the singularity frequencies, the long memory parameter, the WP level and even its location. It is therefore awkward to choose and control in general. To circumvent such a difficulty, a condition involving the singularity frequencies can substitute for the thresholding condition in (\ref{Eq:costthresh1}), that is:
\begin{equation}
\label{Eq:costthresh2}
\EV[W^p_j]=\beta^2_{j,p} \mathds{1}\parenth{\exists~l=1,\ldots,k~|~\nu_l \in I^p_j}
\end{equation}
From the above arguments, it turns out that minimizing the latter cost (with the CW algorithm) will provide us with the same basis as  Algorithm \ref{algo:2}. The main difference is that from a numerical standpoint, our construction algorithm is much faster and stable since there is no need to compute explicitly the band-pass variances, which avoids possible numerical integration problems (because of the PSD singularities).

\section{Analysis of decorrelation properties}
\label{sec:4}
One of the approximations adopted to simulate the Gegenbauer processes using the DWPT is that the coefficients inside each packet of the basis are uncorrelated. Strictly speaking, this is not true, although the expected range of correlation is rather weak as evidently shown by the numerical experiments in Fig.\ref{corrgg}, in contrast to the long-range dependence of the process in the original domain. This section provides a theoretical result that establishes the asymptotic behavior of the covariance between WP coefficients for a $1$-factor Gegenbauer process.
\begin{theorem}
\label{Cov_Wave} If $\psi$ has $q\geq1$ vanishing moments with support $[(N_1-N_2+1)/2,(N_2-N_1+1)/2]$
and $X(t)$ is a stationary $1$-factor Gegenbauer process with Gegenbauer frequency $\nu$. Then the wavelet
packet coefficients covariance $\textrm{Cov}(W_{j_1}^{p_1}(k_1),W_{j_2}^{p_2}(k_2))$ decays as:
\begin{itemize}
  \item $O\left(|2^{j_1}k_1-2^{j_2}k_2|^{2d-1-R_{p_1}-R_{p_2}}\right)$, if $p_1\neq0$ \bf and \rm $p_2\neq0$,
  \item $O\left(|2^{j_1}k_1-2^{j_2}k_2|^{2d-1-R_{\max(p_1,p_2)}}\right)$, if $p_1=0$ \bf or \rm $p_2=0$,
  \item $O\left(|2^{j_1}k_1-2^{j_2}k_2|^{2d-1}\right)$, if $p_1=p_2=0$,
\end{itemize}
for all, $j_1$, $j_2$, $k_1$ and $k_2$ such that $|2^{j_1}k_1-2^{j_2}k_2|> (N^*+1)(2^{j_1}+2^{j_2})$, with
$N^*=\max(N_1,N_2)$, and $R_p=q\sum_{k=0}^{j-1} p^k$ for $p \neq 0$, and $p=\parenth{p^{j-1} p^{j-2} \ldots p^1
p^0}_2$ is the binary representation of $p$. In the last case, we note that $j_1=j_2=j$.
\end{theorem}

Proof:\rm\ {\it See Appendix B}.

This proposition generalizes the results given by \cite{Jen99} and \cite{Jen00} for the case of the FARIMA process. It makes an interesting statement about the order of correlation between well separated WP coefficients, by establishing that the covariance between $W_{j_1}^{p_1}(k_1)$ and $W_{j_2}^{p_2}(k_2)$ decays exponentially over time and scale space. More precisely, the decay speed for $p_1\neq0$ or $p_2\neq0$, depends on the regularity of the wavelet used, on the memory parameter of the process, and indirectly on the location of singularity through the frequency indices $p_1$ and $p_2$. However, keeping the same notations as in Proposition \ref{Cov_Wave}, the larger $q$, the wider the wavelet support and the fewer are the number of wavelet packet coefficients that satisfy the support condition $|2^{j_1}k_1-2^{j_2}k_2|> (N^*+1)(2^{j_1}+2^{j_2})$. Thus, by choosing a wavelet with a large $q$, the rate of decay of autocovariance function increases, but over a subset of WP coefficients. One must then avoid inferring a stronger statement. Nonetheless, the effective support of a wavelet is smaller than the provided bound (see Lemma \ref{lem_support}), and we expect a rapid decay in the WP coefficient's covariance for translations and dilations satisfying $|2^{j_1}k_1-2^{j_2}k_2|> (N^*+1)(2^{j_1}+2^{j_2})$. The following simulation study confirms these remarks.

\section{Simulation results and discussion}
\label{sec:5}

\subsection{Exact correlation of DWPT transformed series}
Suppose we take $\mathbf{X}$ as our input stationary Gegenbauer process vector, whose covariance matrix is $\Gamma$. If $\mathcal{B}$ is the best-ortho-basis provided by our algorithm, it follows that the covariance matrix of the transformed series in the WP domain is:
\begin{equation}
\label{Eq:exactGam}
\Gamma\crochets{\mathcal{B}} = \mathcal{W}^T_{\mathcal{B}} \Gamma \mathcal{W}_{\mathcal{B}}
\end{equation}
where $\mathcal{W}_{\mathcal{B}}$ is the DWPT transform matrix operating on a vector $\mathbf{X}$, whose columns are the basis elements of $\mathcal{B}$. This equation gives the (exact) covariance structure for a given choice of wavelet (type, number of vanishing moments) and treatment of boundaries (e.g. periodic) since both are in $\mathcal{W}_{\mathcal{B}}$.

\begin{figure}[H]
\vspace{-0.25cm}
\hspace{-2.5cm}
\epsfig{file=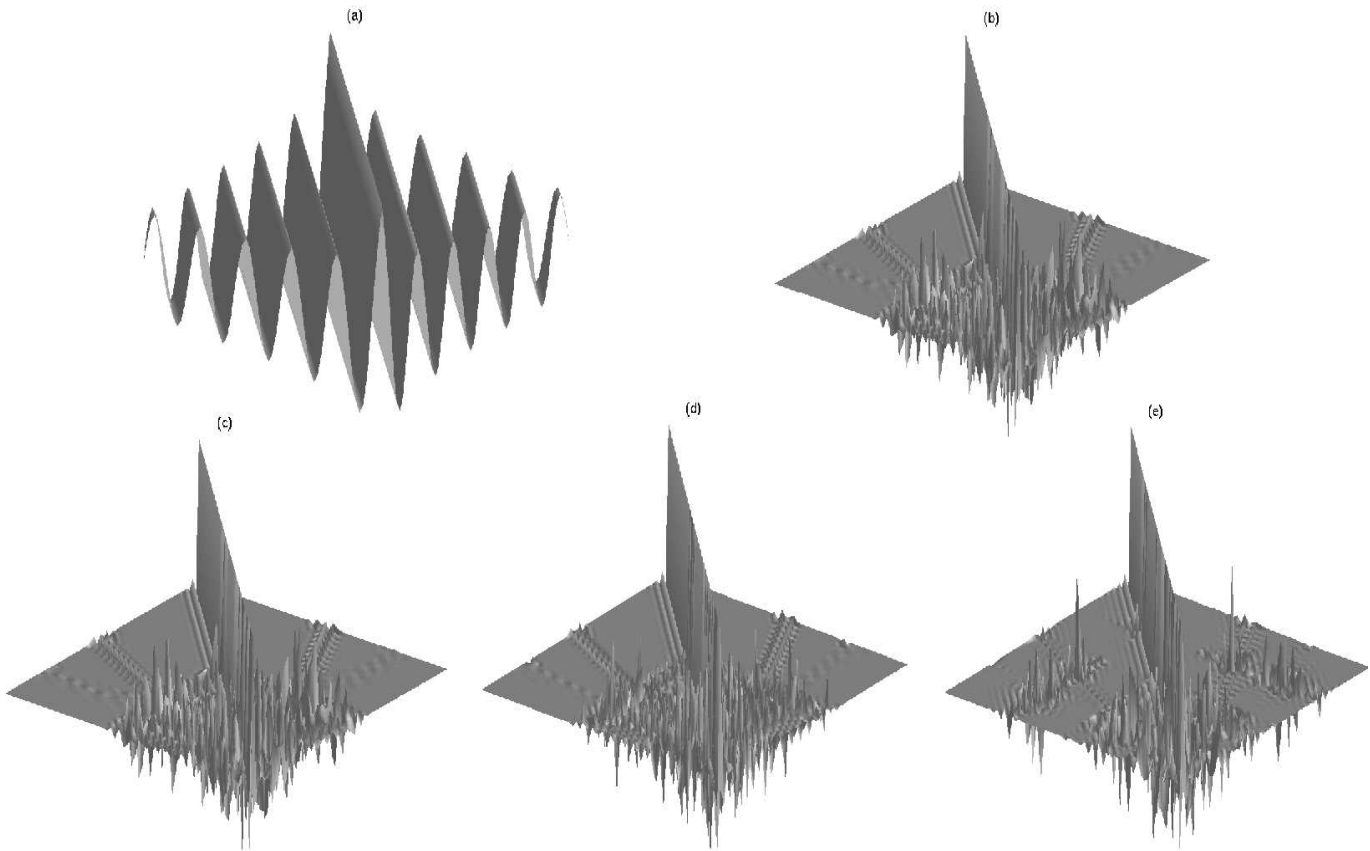,width=1.3\textwidth,height=0.5\textheight}
\caption{(a) Correlation matrix $\Omega$ of a Gegenbauer process with parameters $d=0.4$ and $\nu=1/12$. Correlation matrix $\Omega\crochets{{\mathcal B}}$ of the WP coefficients in our best-ortho-basis for $'db10'$ filter (b), $'sym10'$ filter (c), $'coif5'$ filter (d), and $'bat6'$ filter (e).}
\label{corrgg}
\end{figure}

Fig.\ref{corrgg}.(a) depicts the original correlation matrix $\Omega$ (resulting from $\Gamma$) for a Gegenbauer process vector ($N=64$) with parameters $d=0.4$ and $\nu=1/12$. In Fig.\ref{corrgg}.(b)-(e) are shown the exact correlation matrices resulting from (\ref{Eq:exactGam}), using respectively the Daubechies ($q=10$), Symmlet ($q=10$), Coiflet ($q=10$) and Battle-Lemari\'e wavelets ($q=6$). There is essentially no correlation within the packets that are far from the singularities. The most prominent correlation occurs within the packets close to the singularity. This effect is mainly caused by the support condition stated in Theorem \ref{Cov_Wave} since packets near the singularity are those with smallest length. There is also some correlation between wavelet packets. A significant part of the correlation between two different packets seems to be concentrated along the boundaries between contiguous WP. The latter effect is a consequence of periodic boundary conditions. For example, the periodic boundary effect is higher for the Battle-Lemari\'e spline wavelet, whose support is 59 (compare to the series length of 64). But except boundary effects, this wavelet has a smallest between and intra-packet correlation particularly inside WP close to the singularity. This can be interpreted as a result of a sharper band-pass localization of the Battle-Lemari\'e filters, while the other wavelets have side lobes that yield more energy leaks between bands.

To gain insight into these approximate diagonalizing capabilities of the DWPT, we conduct a larger scale experiment where four Gegenbauer processes with different frequencies and long memory parameters $(d,\nu)$ are studied: three 1-factor with $(0.4,1/12)$, $(0.2,1/12)$, $(0.3,0.016)$ and one 2-factor with $(0.3,1/40)-(0.3,1/5)$. Again the influence of the wavelet on the correlation matrix resulting from (\ref{Eq:exactGam}) is assessed. For comparative purposes, our bases are systematically compared to those of \cite{Whi01} for each process (and wavelet as the best-basis of \cite{Whi01} depends also on the wavelet filter). We also need to consider a criterion to measure the quality of non-correlation. We here propose the Hilbert-Schmidt norm of the diagonalization error, which measures the sum of squares of the off-diagonal elements of the covariance matrix in the best-ortho-basis. As explained above, the method of \cite{Whi01} tends to over-partition the spectral axis yielding to too many packets. Hence, to penalize such configurations and make the comparison fair, we propose the following penalized criterion \cite{DonMalVon98,Donoho97}:
\begin{equation}
\label{Eq:S}
S(\mathcal{B})=\|\Omega\crochets{\mathcal{B}}-\Omega_0\|_{HS}^2+\lambda\#\parenth{\mathcal{B}}
\end{equation}
where $\Omega\crochets{\mathcal{B}}$ is the correlation matrix resulting from (\ref{Eq:exactGam}) and $\Omega_0$ is the correlation matrix of a white noise, i.e. the identity matrix and $\lambda$ is a weight parameter balancing between the diagonalization error and the complexity of the tree associated to $\mathcal{B}$ as measured by $\#(\mathcal{B})$, the number of WP (leaves of the tree) in the basis. 
The value of the weight $\lambda$ is determined by considering two extreme cases. On the one hand, in the Shannon basis, we can assume that the decorrelation of the covariance matrix of the Gegenbauer process is perfect but the tree associated to this basis has too many leaves and the penalty term is the highest; thus $S(\mathcal{B}_S)=2^J\lambda$. On the other hand, if one considers the
basis $\mathcal{B}_0$ composed with only one leaf (i.e. the root packet $W^0_0$), there isn't any decorrelation of the covariance matrix. That is $S(\mathcal{B}_0)=\|\Omega-\Omega_0\|_{HS}^2$, with $\Omega$ the correlation matrix of the Gegenbauer process whose variance-covariance matrix is $\Gamma$. Equating the scores of these two extreme cases yields the following weight:
\begin{equation*}
  \lambda=\frac{\|\Omega-\Omega_0\|_{HS}^2}{2^J-1}.
\end{equation*}

Table.\ref{table:exact} summarizes the scores $S$ obtained for each process as a function of the wavelet filter properties (type and number of vanishing moments). We here assumed time series of length $N=256$. For details about the calculation of the exact autocovariance  function of Gegenbauer processes and hence its associated covariance matrix, see \cite{And86,Chung96,Lapsa97}. These tables show that:
\begin{itemize}
\item The basis provided by our algorithm is systematically better than the one given by \cite{Whi01}, whatever the wavelet and process. Over-partitioning is clearly responsible for the bad performance of the approach in \cite{Whi01}. Meanwhile, the diagonalization error part (not shown here but will be in the next section) remains comparable for both bases. This means that our basis, with a reduced number of packets, does not sacrify the diagonalization quality and yields a diagonalization error comparable to what would be obtained by over-partitioning. It is also worth pointing out that the approach of \cite{Whi01} fails in providing a basis for spline wavelets, and thus cannot be used in this case. The reason is that their best basis search algorithm strongly relies on a threshold on the wavelet packet filter gain, whose choice remains {\it ad hoc}.
\item For a given process, the criterion $S$ decreases as the number of vanishing moments increases. This is in a very good agreement with our expectations as stated in Theorem \ref{Cov_Wave}.
\item From our experiments, we have also noticed that as the number of vanishing moments increases, the best basis provided by \cite{Whi01} tends towards the basis we propose.
\item For all processes, and among all tested wavelets, the Battle-Lemari\'e spline wavelet appears to provide the best score. This confirms our previous observations. Nonetheless, the observed differences between wavelets become less salient at high number of vanishing moments. 
\end{itemize} 

\subsection{Simulation of Gegenbauer processes}
This section is devoted to the illustration of some simulation examples of Gegenbauer processes. The same Gegenbauer processes as in the previous section are considered. For each process, wavelet type and number of vanishing moments, $M=500$ time series of length $N=256$ were generated according to Section \ref{sec:31}, using our basis and that provided by \cite{Whi01} method. For each simulated series, an unbiased estimate of the autocovariance function for the first $N/2$ lags was calculated. An average of the autocovariance function (over the $M$ estimates) was then obtained and the associated correlation matrix $\bar{\Omega}$ was constructed. Finally, the HS norm of deviation between the true and averaged sample correlation matrices was computed:
\[
B(\mathcal{B})=\|\Omega-\bar{\Omega}\|_{HS}^2
\]
As previously, a penalized version of $B(\mathcal{B})$ by the complexity of the tree associated to $\mathcal{B}$, as in (\ref{Eq:S}) was also calculated (denoted $B_{\text{pen}}$)\footnote{Note that for our best-basis algorithm, and for a given process, the scores $B_{\text{pen}}$ is simply $B$ plus a constant for all wavelets, as the penalty part in $B_{\text{pen}}$ only depends on the singularity frequencies.}. In order to determine which part of the score $B_{\text{pen}}$ is the largest contributor to the performance, and in order to not favour our best basis construction algorithm, both $B$ and $B_{\text{pen}}$ are displayed. The score $B$ of the Hosking method \cite{Hosk84}, which is an exact simulation scheme, is also reported. The results are summarized in Table.\ref{table:simu}.

As revealed by these tables, the deviation error part $B$ is comparable between the two best basis construction methods, but the penalized version differs significantly. This is caused by a fairly large difference in the "size" of the basis. Again, this backs up the statement that the method of \cite{Whi01} over-partitions the spectrum, and also agrees with the fact that in terms of performance our method generates as reasonable Gegenbauer processes as \cite{Whi01} with less number of packets. This also clearly provides a numerical support to our claim that good quality DWPT-based best-basis search, and then simulation, of Gegenbauer processes can be achieved without necessarily depending on the wavelet choice, just as it has been extensively done for $1/f$ processes using the DWT. But, one has to keep in mind that the quality of the reconstructed covariance structure (by assuming almost decorrelation of WP coefficients in the best-ortho-basis), compared to the true covariance of a Gegenbauer process will still depend on the wavelet. From this point of view (decorrelation performance), the numerical results observed for simulated data essentially confirm those reported in the previous subsection.

Both the score $B$ and its penalized version exhibit a decreasing tendency with increasing number of vanishing moments. This numerical evidence is a confirmation of the previous subsection findings and support our claims in Theorem \ref{Cov_Wave}. The Battle-Lemari\'e spline wavelet seems to perform the best (in terms of both $B$ and $B_{\text{pen}}$), followed closely by the symmlets. The difference in performance between all wavelet types vanishes as $q$ increases.

\section{Conclusion}
In this paper, we provided a new method to build approximate diagonalizing bases for $k$-factor Gegenbauer
processes. Exploiting the intuitive fact that a wavelet packet library contains the basis where a Gegenbauer
process could be (almost) whitened, our best-ortho-basis search algorithm was formulated in the case of
$1$-factor process and the fast search algorithm of Coifman-Wickerhauser was adapted to find this best basis.
Using this framework, our methodology was posed in a well principled way and the uniqueness of the basis was
guaranteed. Furthermore, unlike the approach \cite{Whi01}, it is very fast (see simulations), does not depend on the
wavelet choice, and is not very sensitive to the length of the time series. As the method construction of the
best basis for simulation of a $k$-factor Gegenbauer process relies on the $1$-factor construction method, the
same conclusions hold. 

Then, we studied the error of diagonalization in the best-ortho-basis. Towards this goal, we established the
decay speed of the correlation between two WP coefficients. These results generalize the work of
\cite{Jen99} and \cite{Jen00} provided in the case of FARIMA processes. The numerical evidence shown by our experimental study 
confirmed these theoretical findings. It has also shown that the algorithm introduced in the paper is appealing in that it provides good quality simulated Gegenbauer processes with computational simplicity and reduced complexity bases independently of the wavelet, which is a clear improvement over the existing method in \cite{Whi01}. Owing to these appealing theoretical and empirical properties, and given its practical simplicity, we feel the general practitioner will be attracted to our simulator.

This new method of simulating Gegenbauer processes gives a new perspective for analyzing processes whose
PSD singularities occur at any frequency in the Nyquist interval. In such a task, one could have
the basis by knowing the process parameters ($\nu$ in particular). Thus, our method has a direct application for bootstrap-based inference in the presence of Gegenbauer noise.

A remaining important open problem is how we could extend this work if the
question of interest becomes that of estimating the parameters of a $k$-factor Gegenbauer given one or more sample paths
of this process. This estimation problem can be accomplished in a maximum likelihood framework once
the diagonalizing basis is found. In this case, the best-ortho-basis cannot be found by a naive straightforward
application of Algorithm \ref{algo:2}. Nevertheless, we have some promising directions that
are now under investigation. 
Establishing the asymptotic behavior of such estimators also remains an open problem. One could also refine the estimation process by handling the residual correlation structure of the WP coefficients via explicit modeling by a low-order autoregressive process as recently suggested in \cite{Craigmile04} for $1/f$ fractionally-differenced processes. Additional research is still required and our current work is focusing on these directions.

\section*{Appendix A}
\small
\proof{Proposition \ref{Prop_1_fact}}
\begin{itemize}
\item Let's consider the node $(j,p)$. We compute the variance of the WP coefficients at its two
  children: $(j+1,2p)$ and $(j+1,2p+1)$. Without loss of generality, we assume that the frequency $\nu$ is in the interval
  $I_{j+1}^{2p}=[\frac{2p}{2^{j+1}},\frac{2p+1}{2^{j+1}}[$. Then a good approximation of the variance of the
  WP coefficient is given by the integral over the interval $I_{j+1}^{2p}$ of the PSD. On
  this interval, a very good approximation to the PSD of the process
  $f(\lambda)=\frac{\sigma^2}{2\pi}|2(\cos2\pi\lambda-\cos2\pi\nu)|^{-2d}$ is given by $C_0|\lambda-\nu|^{-2d}$
  with $C_0$ a positive constant.

Two different cases are then distinguished with associated values of $A_0$:

$\ast$ \underline{Case $\nu\leq\frac{4p+1}{2^{j+2}}$:}
\begin{eqnarray*}
    \mathbb{V}[W_{j+1}^{2p}] &=& C_0\int_{\frac{2p}{2^{j+1}}}^{\nu}|\nu-\lambda|^{-2d}d\lambda+C_0\int_{\nu}^{\frac{2p+1}{2^{j+1}}}|\lambda-\nu|^{-2d}d\lambda\\   &=&\frac{C_0}{1-2d}\left(\left(\nu-\frac{2p}{2^{j+1}}\right)\left(\left(\nu-\frac{2p}{2^{j+1}}\right)^2\right)^{-d}+\left(\frac{2p+1}{2^{j+1}}-\nu\right)\left(\left(\frac{2p+1}{2^{j+1}}-\nu\right)^2\right)^{-d}\right)\\
 &=&\frac{C_0}{1-2d} u^{1-2d}\left(1-\left(1-\frac{1}{2^{j+1} u}\right)^{1-2d}\right), \ \ \text{where~} u=\frac{2p+1}{2^{j+1}}-\nu
\end{eqnarray*}
\begin{eqnarray*}
    &\geq&\frac{C_0}{2^{j+1}}u^{-2d}\left(1+\frac{2d}{2^{j+1}u}\right)\\ &\geq&\frac{C_0}{2^{j+1}}u^{-2d}\left(1+\frac{2d}{2^{j+1}\frac{1}{2^{j+1}}}\right)=\frac{C_0}{2^{j+1}}\left(\frac{2p+1}{2^{j+1}}-\nu\right)^{-2d}\left(1+2d\right)\\
\end{eqnarray*}

where the last inequality is a consequence of the fact that $u \leq 2^{-(j+1)}$ in this case.

To compute the variance of $W_{j+1}^{2p+1}$ we denote $\lambda^*$ the location of the maxima of the PSD $f$ over the interval $I_{j+1}^{2p+1}$. As $f$ is a non-increasing function over $[\frac{2p+1}{2^j},\frac{2p+2}{2^{j+1}}]$, it follows that this variance is bounded by a rectangle area (to a good approximation in this case):
\begin{eqnarray*}
  \mathbb{V}[W_{j+1}^{2p+1}]
  &\leq& \frac{\sigma^2}{2\pi2^{j+1}}|2(\cos2\pi\lambda^*-\cos2\pi\nu)|^{-2d}.
\end{eqnarray*}
Using the same approximation of the PSD as previously, we obtain:
$$\mathbb{V}[W_{j+1}^{2p+1}]\leq\frac{C_0}{2^{j+1}}\left(\frac{2p+1}{2^{j+1}}-\nu\right)^{-2d}.$$
Thus
$$\mathbb{V}[W_{j+1}^{2p+1}]\leq\frac{1}{1+2d}\mathbb{V}[W_{j+1}^{2p}]=A_0\mathbb{V}[W_{j+1}^{2p}],~ 0<A_0<1$$
Therefore, in this case we can write that $\mathbb{V}[W_{j+1}^{2p+1}]\ll\mathbb{V}[W_{j+1}^{2p}]$, and following the criterion defined in section \ref{sec:31} we have $\EV[W_{j+1}^{2p+1}]=0$. Consequently, at the node $(j,p)$, the algorithm (\ref{new_algo_CW}) gives us:
\begin{equation}
\mathcal{B}_j^p=\mathcal{B}_{j+1}^{2p}\cup\mathcal{B}_{j+1}^{2p+1}.
\end{equation}

$\ast$ \underline{Case $\nu\geq\frac{4p+1}{2^{j+2}}$:}
Using the same steps as in the first case, we prove that:
\begin{eqnarray}\label{var_gauche}
    \mathbb{V}[W_{j+1}^{2p}] &=&
    \frac{C_0}{1-2d}\left(\frac{2p+1}{2^{j+1}}-\nu\right)^{1-2d}\left(1-\left(1-\frac{1}{2^{j+1}(\frac{2p+1}{2^{j+1}}-\nu)}\right)^{1-2d}\right).
  \end{eqnarray}
When $\frac{2p+1}{2^{j+1}}\sim\nu$, the rectangle approximation is no longer valid. But a good approximation of the PSD in the interval $[\frac{2p+1}{2^{j+1}},\frac{2p+2}{2^{j+1}}]$ can be $C_0|\lambda-\nu|$, where the constant $C_0$ is the same as previously. Then, after some manipulations:
\begin{eqnarray}
\label{var_droite}
    \mathbb{V}[W_{j+1}^{2p+1}] &=&
    \frac{C_0}{1-2d}\left(\frac{2p+1}{2^{j+1}}-\nu\right)^{1-2d}\left(1+\left(1-\frac{1}{2^{j+1}\left(\frac{2p+1}{2^{j+1}}-\nu\right)}\right)^{1-2d}\right).
  \end{eqnarray}
  
As by assumption $\frac{2p+1}{2^{j+1}}\sim\nu$, we have that $0<\frac{2p+1}{2^{j+1}}-\nu<\frac{1}{2^{j+2}}$ and
then,
$$\left(1-\frac{1}{2^{j+1}\left(\frac{2p+1}{2^{j+1}}-\nu\right)}\right)^{1-2d}<0.$$
Finally, combining equations (\ref{var_gauche}) and (\ref{var_droite}), we obtain $\mathbb{V}[W_{j+1}^{2p+1}] \sim A_1\mathbb{V}[W_{j+1}^{2p}],$ where:
\begin{equation}
    A_1 =   \frac{1+\left(1-\frac{1}{2^{j+1}\left(\frac{2p+1}{2^{j+1}}-\nu\right)}\right)^{1-2d}}{1-\left(1-\frac{1}{2^{j+1}\left(\frac{2p+1}{2^{j+1}}-\nu\right)}\right)^{1-2d}} < 1.
\end{equation}
In this case we can write that $\mathbb{V}[W_{j+1}^{2p+1}]\ll\mathbb{V}[W_{j+1}^{2p}]$, and using the criterion defined in section \ref{sec:31} we obtain $\EV[W_{j+1}^{2p+1}]=0$. Finally, at the node $(j,p)$, algorithm (\ref{new_algo_CW}) gives us:
\begin{equation}
\mathcal{B}_j^p=\mathcal{B}_{j+1}^{2p}\cup\mathcal{B}_{j+1}^{2p+1}.
\end{equation}

\item In the case where the frequency $\nu$ is in the closure of the intervals $I_{j+1}^{2p}$ and
$I_{j+1}^{2p+1}$, we have no relationship as $\EV[W_{j+1}^{2p}]=0$ or $\EV[W_{j+1}^{2p+1}]=0$, 
and one could not conclude. However, fortunately, at the depth
$j+2$, we still have:
\[
\EV[W_{j+2}^{4p+1}]=0\ \ \ \ \ \textrm{and}\ \ \ \ \
\EV[W_{j+2}^{4p+2}]=0.
\]
Then we easily obtain that for algorithm (\ref{new_algo_CW}):
$$\mathcal{B}_j^p=\mathcal{B}_{j+2}^{4p}\cup\mathcal{B}_{j+2}^{4p+1}\cup\mathcal{B}_{j+2}^{4p+2}\cup\mathcal{B}_{j+2}^{4p+3}.$$
\end{itemize}
\endproof

\proof{Proposition \ref{Prop_k_fact}}
\begin{enumerate}
  \item Let $(j,p)$ be a node. We assume that this node is not in the tree $\mathcal{T}'$.
  It means that this node is not in the tree $\mathcal{T}_1$ neither in
  $\mathcal{T}_2$ and in terms of threshold, we have
  $$\EV[W_{j+1}^{2p}(1)]=0\ \ \textrm{or}\ \
    \EV[W_{j+1}^{2p+1}(1)]=0\ \ \ \ \ \ \
    \ \ \textrm{\bf and\rm}\ \ \ \ \ \ \ \ \
    \EV[W_{j+1}^{2p}(2)]=0\ \ \textrm{or}\ \
    \EV[W_{j+1}^{2p+1}(2)]=0.$$
  As the tree $\mathcal{T}_2$ is associated to the best-ortho-basis of a $1$-factor Gegenbauer process, the fact
  that the node $(j,p)$ is node in the tree $\mathcal{T}_2$ means that the frequency $\nu_k$ is not in the
  interval $I_{j}^p=[\frac{p}{2^j},\frac{p+1}{2^j}]$. Then in this interval, the function
  $|2(\cos2\pi\lambda-\cos2\pi\nu_k)|^{-2d_k}$ is bounded and has a maximum at frequency $\lambda^*\in
  I_j^p$. Then,
  \begin{eqnarray}
  \label{Eq:varkfactor}
    \mathbb{V}[W_{j+1}^{2p}]
    &\leq&
    \frac{\sigma^2}{2\pi}|2(\cos2\pi\lambda^*-\cos2\pi\nu_k)|^{-2d_k}\int_{\frac{2p}{2^{j+1}}}^{\frac{2p+1}{2^{j+1}}}\prod_{i=1}^{k-1}|2(\cos2\pi\lambda-\cos2\pi\nu_i)|^{-2d_i}d\lambda \nonumber \\
    &=&
    \frac{\sigma^2}{2\pi}|2(\cos2\pi\lambda^*-\cos2\pi\nu_k)|^{-2d_k}\mathbb{V}[W_{j+1}^{2p}(1)].
  \end{eqnarray}
  and, using the same argument,
  $$\mathbb{V}[W_{j+1}^{2p+1}]\leq\frac{\sigma^2}{2\pi}|2(\cos2\pi\lambda^*-\cos2\pi\nu_k)|^{-2d_k}\mathbb{V}[W_{j+1}^{2p+1}(1)].$$
  Finally, as $\EV[W_{j+1}^{2p}(1)]=0$ or $\EV[W_{j+1}^{2p+1}(1)]=0$, we
  have $\mathbb{V}[W_{j+1}^{2p}]\gg\mathbb{V}[W_{j+1}^{2p+1}]$ or $\mathbb{V}[W_{j+1}^{2p}]\gg\mathbb{V}[W_{j+1}^{2^rp+1}],$
  which means,
  $$\EV[W_{j+1}^{2p}]=0\ \ \ \ \ \ \ \textrm{or}\ \ \ \ \ \ \EV[W_{j+1}^{2p+1}]=0.$$
  Then, $\mathcal{B}_j^p=\mathcal{B}_{j+1}^{2p}\cup\mathcal{B}_{j+1}^{2p+1}$ and so the node $(j,p)$ is not in the tree $\mathcal{T}$. Finally,
  $$\mathcal{B}\subset\mathcal{B}'.$$


  \item Here $(j,p)$ and $(j+r^*,2^{r^*}p+s^*)$ (for $s^*=0,\dots,2^{r^*}-1$) are in the tree $\mathcal{T}'$. We
  denote $r$ the minimum value of $r^*$
  for which there exists a $s$ ($s=0,\dots,2^r-1$) such that the node $(j+r,2^{r}p+s)$ is in the tree
  $\mathcal{T}'$.\\
  Then the fact that the nodes $(j,p)$ and $(j+r,2^rp+s)$ are in the tree $\mathcal{T}'$ means that $(j,p)$ is in
  $\mathcal{T}_1$ or in $\mathcal{T}_2$ and $(j+r,2^rp+s)$ is in $\mathcal{T}_2$ or in $\mathcal{T}_1$ (it is
  important to remark that both $(j,p)$ and $(j+r,2^rp+s)$ cannot be in $\mathcal{T}_1$ or in
  $\mathcal{T}_2$). Without loss of generality, we assume that $(j,p)$ is in $\mathcal{T}_2$ and
  $(j+r,2^*p+s)$ is in $\mathcal{T}_1$. All the calculations made in the following remain valid if
  we consider that $(j,p)$ is in $\mathcal{T}_1$ and $(j+r,2^rp+s)$ is in $\mathcal{T}_2$. To
  simplify the notations, we assume also that there exists a $s$ which is even.\\
  We denote $W_j^p(1)$ and $W_j^p(2)$, for $j=0,\dots,J$ and $p=0,2^j-1$, the wavelet packet coefficients of
  respectively the processes $(X^1_t)_t$ and $(X^2_t)_t$. From these sub-processes, we have that for the
  algorithm CW,
  \begin{itemize}
    \item for the tree $\mathcal{T}_1$:
    $$\mathcal{B}_{j+r-1}^{2^{r-1}p+\frac{s}{2}}(1)=\mathcal{B}_{j+r}^{2^rp+s}(1)\cup\mathcal{B}_{j+r}^{2^rp+s+1}(1)$$
    because $\EV[W_{j+r}^{2^rp+s}(1)]=0$ or
    $\EV[W_{j+r}^{2^rp+s+1}(1)]=0$,
    \item for the tree $\mathcal{T}_2$: $$\mathcal{B}_j^p(2)=\mathcal{B}_j^p(2)$$
    because 
    $\EV[W_{j+1}^{2p}(2)]=\mathbb{V}[W_{j+1}^{2p}(2)]$ and
    $\EV[W_{j+1}^{2p+1}(2)]=\mathbb{V}[W_{j+1}^{2p+1}(2)]$,
  \end{itemize}
  We consider the intervals $I_{j+r}^{2^rp+s}=[\frac{2^rp+s}{2^{j+r}}, \frac{2^rp+s+1}{2^{j+r}}]$ and
  $I_{j+r}^{2^rp+s+1}=[\frac{2^rp+s+1}{2^{j+r}}, \frac{2^rp+s+2}{2^{j+r}}]$. As the node $(j,p)$ is in the tree
  $\mathcal{T}_2$, and as $\mathcal{T}_2$ the tree of a basis, the frequency $\nu_k$ is not in the interval
  $I_{j+r}^{2^rp+s}\cup I_{j+r}^{2^rp+s+1}$.\\
  We denote $\lambda^*$ the location of the maximum of $|2(\cos2\pi\lambda-\cos2\pi\nu_k)|^{-2d_k}$ in the
  interval $I_{j+r-1}^{2^{r-1}p+s/2}$ (Note that the maximum is bounded).
  From (\ref{Eq:varkfactor}), we have:
  \begin{eqnarray*}
    \mathbb{V}[W_{j+r}^{2^rp+s}]
    &\leq&
    \frac{\sigma^2|2(\cos2\pi\lambda^*-\cos2\pi\nu_k)|^{-2d_k}}{2\pi}\mathbb{V}[W_{j+r}^{2^rp+s}(1)],
  \end{eqnarray*}
  \begin{eqnarray*}
    \mathbb{V}[W_{j+r}^{2^rp+s+1}]
    &\leq&
    \frac{\sigma^2|2(\cos2\pi\lambda^*-\cos2\pi\nu_k)|^{-2d_k}}{2\pi}\mathbb{V}[W_{j+r}^{2^rp+s+1}(1)].
  \end{eqnarray*}
  Then, as $\EV[W_{j+r}^{2^rp+s}(1)]=0$ or $\EV[W_{j+r}^{2^rp+s+1}(1)]=0$,
  we have $\mathbb{V}[W_{j+r}^{2^rp+s}]\gg\mathbb{V}[W_{j+r}^{2^rp+s+1}]$ or $\mathbb{V}[W_{j+r}^{2^rp+s}]\gg\mathbb{V}[W_{j+r}^{2^rp+s+1}],$
  which means,
  $$\EV[W_{j+r}^{2^rp+s}]=0\ \ \ \ \ \ \ \textrm{or}\ \ \ \ \ \ \EV[W_{j+r}^{2^rp+s+1}]=0.$$
  Then, because of the particular choice of $r$, 
  $$\mathcal{B}_j^p=\bigcup_{i=0}^{2^r-1}\mathcal{B}_{j+r}^{2^rp+i}.$$
  Finally, the node $(j,p)$ is not in the tree $\mathcal{T}$.\\
  However, the fact that the node $(j+r,2^rp+s)$ is in the tree $\mathcal{T}_1$ means that,
  $$\EV[W_{j+r+1}^{2^{r+1}p+2s}(1)]=W_{j+r+1}^{2^{r+1}p+2s}(1)\ \ \ \ \textrm{and}\ \ \ \ \ \EV[W_{j+r+1}^{2^{r+1}p+2s+1}(1)]=W_{j+r+1}^{2^{r+1}p+2s+1}(1).$$
  As the frequency $\nu_k$ is not in the interval $I_{j+r}^{2^rp+s}$, we obtain easily that
  $$\EV[W_{j+r+1}^{2^{r+1}p+2s}]=W_{j+r+1}^{2^{r+1}p+2s}\ \ \ \ \textrm{and}\ \ \ \ \ \EV[W_{j+r+1}^{2^{r+1}p+2s+1}]=W_{j+r+1}^{2^{r+1}p+2s+1}.$$
  Finally, the node $(j+r,2^rp+s)$ is in the tree $\mathcal{T}$. Using the argument, we show that the node
  $(j+r^*,2^{r^*}p+s^*)$ is in the tree $\mathcal{T}$.
\end{enumerate}
\endproof

\vspace{-1cm}
\section*{Appendix B}
To prove Theorem \ref{Cov_Wave}, the following preliminary lemmas are needed.
\begin{lemme}
  \label{lem_moment}
  Let $\psi$ be a wavelet with $q$ vanishing moments, and the associated high-pass QMF filters $h$ and $g$.
  Then, for all $j$ and $p=0,\dots,2^j-1$, the moments of the WP function $\psi^j_p$ are such that:
  \begin{equation*}
  \mathcal{M}_{j,p}(r)=\int_{\mathbb{R}} t^r \psi_j^p(t) dt = \delta(r)\delta(p), ~ for ~ 0 \leq r < R_p
  \end{equation*}
  where $R_0=1$ and $R_p=q\sum_{k=0}^{j-1} p_k$ for $p \neq 0$, and $p=\parenth{p_{j-1} p_{j-2} \ldots p_1 p_0}_2$ is the binary representation of $p$.
\end{lemme}

\proof{Lemma \ref{lem_moment}}
Note that for $p=1$ (wavelet basis), our result specializes to the traditional relation $\mathcal{M}_{j,1}(r)=0$
for $0 \leq r < q$. The lemma can be proved either by induction in the original domain, or using explicit proof
in the Fourier domain. We shall proceed according to the latter.
By iterating the actions of the QMF filters $h$ or $g$, from
the root of the binary tree, to extract the appropriate range of frequencies, one can write that:
\begin{equation}
\label{Eq:psijp1}
\hat{\psi}_{j}^p(\omega) = \crochets{\prod_{k=0}^{j-1}
F_{p_{j-k-1}}(2^{k}\omega)} \hat{\phi}(\omega)
\end{equation}
where the sequence of filters $F_{p_k}$ is chosen according to $p=2^{j-1}p_{j-1} + 2^{j-2} p_{j-2}+ \ldots + 2p_1 + p_0$:
\begin{equation}
F_{p_k} = \begin{cases}
    \hat{h} & if ~ p_k = 0 \\
    \hat{g} & if ~ p_k = 1
      \end{cases}
\end{equation}
and $\hat{\phi}(0) \neq 0$.

For compactly supported wavelets with $q$ vanishing moments, the associated high-pass filters $\hat{g}$
has $q-1$ zeros at $\omega=0$:
\begin{equation}
\label{Eqqmfg} 
\hat{g}(\omega) = \parenth{1-e^{-i\omega}}^q P\parenth{e^{i\omega}}
\end{equation}
where $P(.)$ is a trigonometric polynomial bounded around $\omega=0$. The number of vanishing moments of
$\psi_j^p(t)$ is equivalently given by the number of vanishing derivatives of $\hat{\psi}_j^p(\omega)$ at
$\omega=0$, that is:
\begin{equation}
\mathcal{M}_{j,p}(r)=\crochets{\parenth{\frac{1}{i}\partial_\omega}^r \hat{\psi}_j^p(\omega)}_{\omega=0}  ~
\text{for} ~ r=0,\ldots,R_p-1
\end{equation}

\begin{itemize}
\item If $p=0$, $\hat{\psi}_{j}^p(\omega)$ is just the product of low-pass filters, and
$\psi_j^0(t)=\phi_j(t)$ the scaling function at depth $j$. Then, $\mathcal{M}_{j,0}(r)=\hat{\phi}_j(0)$, which
is non-zero with $R_0=1$. If additional constraints are imposed on the wavelet choice
(e.g. Coiflets), $\mathcal{M}_{j,0}(r)$ might be zero for $1 \leq r < q$.

\item If $p\neq0$, from (\ref{Eq:psijp1}) we can write:
\begin{equation}
\hat{\psi}_{j}^p(\omega) = \prod_{k|p_{j-k-1}=1} \parenth{1-e^{-i2^{k}\omega}}^q Q(\omega)
\end{equation}
where $Q(.)$ is again bounded around $\omega=0$. The number of vanishing moments is then given by the number of
zeros at $\omega=0$ which is $R_p=q\sum_k p_k$. The lemma follows.
\end{itemize}
\endproof

\begin{lemme}
  \label{lem_support}
If the QMF $h$ has a support in $\crochets{N_1,N_2}$, then the support of the WP function
$\psi_j^p(t)$ at each node $\parenth{j,p}$ in the WP binary tree is always included in
$\crochets{-2^j\parenth{N^*+1},2^j\parenth{N^*+1}}$, with $N^*=\max\parenth{|N_1|,|N_2|}$.
\end{lemme}

\proof{Lemma \ref{lem_support}}
This is proved by induction. We also use the fact that $\psi_0$ will be supported in the interval
$\crochets{N_1,N_2}$ and $\psi_1$ in $\crochets{\frac{N_1-N_2+1}{2},\frac{N_2-N_1+1}{2}}$ (see e.g. Mallat (1998)\cite{Mal98}, Proposition 7.2).
\endproof

\begin{lemme}
\label{mom_acf} Let $\mathcal{I}$ be a collection of disjoint dyadic intervals $I^j_p$ whose union is the
positive half line, and $\mathcal{B}=\{\psi_j^p(t-2^jk):0 \leq k < 2^{J-j}, I_j^p \in \mathcal{I}\}$ is the
associated orthonormal basis. Let $h$ and $g$ the QMFs as defined in (\ref{Eqqmfg}). The vanishing moments
$\mathcal{M}_{j_1,j_2}^{p_1,p_2}(m)$ of the inter-correlation function $\Lambda_{j_1,j_2}^{p_1,p_2}(h)$ of
$\psi^{p_1}_{j_1}(t)$ and $\psi^{p_2}_{j_2}(t) \in \mathcal{B}$ satisfy:

\begin{tabular}{llll}
1) $p_1 \neq 0$ {\bf and} $p_2 \neq 0$: & $\mathcal{M}_{j_1,j_2}^{p_1,p_2}(m) = 0$ & for &  $0 \leq m<R_{p_1}+R_{p_2}$, \\
2) $p_1 \neq 0$ {\bf or}  $p_2 \neq 0$:  & $\mathcal{M}_{j_1,j_2}^{p_1,p_2}(m) = 0$ & for &  $0 \leq m< R_{\max(p_1,p_2)}$, \\
3) $p_1=p_2=0$: & $\mathcal{M}_{j_1,j_2}^{p_1,p_2}(m) = 0$ & for & $1 \leq m<2q$.
\end{tabular}
\\
Furthermore, the support of $\Lambda_{j_1,j_2}^{p_1,p_2}(h)$ is included in
$\crochets{-\parenth{N^*+1}\parenth{2^{j_1}+2^{j_2}},\parenth{N^*+1}\parenth{2^{j_1}+2^{j_2}}}$.
\end{lemme}

\proof{Lemma \ref{mom_acf}}
By definition of the inter-correlation function, we have:
\begin{equation}
\Lambda_{j_1,j_2}^{p_1,p_2}(h)=\int \psi_{j_1}^{p_1}(t)\psi_{j_2}^{p_2}(t-h)dt
\end{equation}
As these WP functions belong to the orthonormal basis $\mathcal{B}$, then at integer lags $\Lambda_{j_1,j_2}^{p_1,p_2}(n)=\delta\parenth{j_1-j_2}\delta\parenth{p_1-p_2}\delta\parenth{n}$.

As far as the support is concerned, it is not a difficult matter to see, using Lemma \ref{lem_support}, that
$\Lambda_{j_1,j_2}^{p_1,p_2}(h)$ is supported in
$\crochets{-\parenth{N^*+1}\parenth{2^{j_1}+2^{j_2}},\parenth{N^*+1}\parenth{2^{j_1}+2^{j_2}}}$.

Let's now turn to the moments of $\Lambda_{j_1,j_2}^{p_1,p_2}(h)$.
\begin{enumerate}
\item $p_1 \neq 0$ {\bf and} $p_2 \neq 0$: 
\label{sit1}\\
In this case, we know that $\psi_{j_1}^{p_1}$ and $\psi_{j_2}^{p_2}$ have
respectively $R_{p_1}$ and $R_{p_2}$ vanishing moments as defined in Lemma \ref{lem_moment}. Then,
$\Lambda_{j_1,j_2}^{p_1,p_2}(h)$ will have $R_{p_1}+R_{p_2}$ vanishing moments since,
\begin{align}
  \int h^m\Lambda_{j_1,j_2}^{p_1,p_2}(h)dh &= -\int\int(v-u)^m\psi_{j_1}^{p_1}(v)\psi_{j_2}^{p_2}(u)dudv \nonumber\\
  &= -\sum_{n=0}^m \parenth{-1}^n {m\choose n} \int v^{m-n} \psi_{j_1}^{p_1}(v)dv \int u^n\psi_{j_2}^{p_2}(u)du = 0,
\end{align}
for $0 \leq m<R_{p_1}+R_{p_2}$. Here, we used uniform convergence and continuity to invert the order of summation and integration. Note that the Fubini theorem allows us to invert the order of integrals.
\item $p_1 \neq 0$ {\bf or} $p_2 \neq 0$:
\label{sit2}
Without loss of generality, assume that $p_1 \neq 0$ and $p_2=0$. The same reasoning as above can be adopted to
conclude that for $0 \leq m<R_{p_1}$:
\begin{align}
\label{Eqsit2}
   \int h^m\Lambda_{j_1,j_2}^{p_1,0}(h)dh &= -\sum_{n=0}^m \parenth{-1}^n {m\choose n} \int v^{m-n} \psi_{j_1}^{p_1}(v)dv \int u^n\phi_{j_2}(u)du = 0.
\end{align}
\item $p_1=p_2=0$: 
\label{sit3}\\
In an orthonormal basis of wavelet packets, this situation is not possible unless $j_1 = j_2 = j$. Thus,
\begin{align}
  \int h^m\Lambda_{j,j}^{0,0}(h)dh &= 2^{-j}\int\int h^m\phi\parenth{2^{-j}t}\phi\parenth{2^{-j}\parenth{t-h}}dtdh \nonumber\\
  &= 2^{j(m+1)}\int\int u^m\phi(v)\phi(v-u)dudv = 0,
\end{align}
for $1 \leq m<2q$, where the latter result is proved in \cite{Bey92}.
\end{enumerate}
\endproof

\proof{Theorem \ref{Cov_Wave}}
Here we are interested in the covariance between the WP coefficients $W_{j_1}^{p_1}(k_1)$ and
$W_{j_2}^{p_2}(k_2)$. We have:
\begin{eqnarray}
  \textrm{Cov}\left[W_{j_1}^{p_1}(k_1),W_{j_2}^{p_2}(k_2)\right]
  &=& \int\int\mathbb{E}[X(t)X(s)]\psi_{j_1}^{p_1}(t-2^{j_1}k_1)\psi_{j_2}^{p_2}(s-2^{j_2}k_2)dtds \nonumber\\
  &=& \int\int\cos\left(\nu(t-s)\right)|t-s|^{2d-1}\psi_{j_1}^{p_1}(t-2^{j_1}k_1)\psi_{j_2}^{p_2}(s-2^{j_2}k_2)dtds
\end{eqnarray}
After three changes of variables, $u=t-2^{j_1}k_1$ and $v=s-2^{j_2}k_2$, then $u=t'$ and $v=t'-h$ and finally $\alpha=2^{j_1}(k_1-2^{j_2-j_1}k_2)$, we obtain:
\begin{equation}
\textrm{Cov}\left[W_{j_1}^{p_1}(k_1),W_{j_2}^{p_2}(k_2)\right]=\int\cos(\nu(h+\alpha))|h+\alpha|^{2d-1}\Lambda_{j_1,j_2}^{p_1,p_2}(h)dh.
\end{equation}
From Lemma \ref{lem_support} we know that that the support of $\Lambda_{j_1,j_2}^{p_1,p_2}(h)$ is included in $[-(2^{j_1}+2^{j_2})(N^*+1),(2^{j_1}+2^{j_2})(N^*+1)]$. As $h$ is in the support of $\Lambda_{j_1j_2}^{p_1p_2}$ and by assumption $\alpha> (N^*+1)(2^{j_1}+2^{j_2})$, we have $h/\alpha<1$. Hence, using the binomial series expansion of $\left|1+\frac{h}{\alpha}\right|^{2d-1}$ and the fact that $\cos(\nu(\alpha+h)) \sim \cos(\nu\alpha)$ for large $\alpha$, it follows that:
\begin{equation}
\label{eq_general} \textrm{Cov}\left[W_{j_1}^{p_1}(k_1),W_{j_2}^{p_2}(k_2)\right] \sim
|\alpha|^{2d-1}\cos(\nu\alpha)\left\{\int\Lambda_{j_1,j_2}^{p_1,p_2}(h)dh+\sum_{i=1}^{\infty}{2d-1\choose
i}\int\left(\frac{h}{\alpha}\right)^i\Lambda_{j_1,j_2}^{p_1,p_2}(h)dh\right\}.
\end{equation}
We must then provide an upper bound on the integrals inside the braces. In the following we distinguish three
different cases depending on the number of vanishing moments of $\Lambda_{j_1,j_2}^{p_1,p_2}$ according to Lemma
\ref{mom_acf}, that is:
\begin{enumerate}
\item If $p_1 \neq 0$ {\bf and} $p_2 \neq 0$, then $\mathcal{M}^{j_1,j_2}_{p_1,p_2}(m) = 0$, for $0\leq
m<R_{p_1}+R_{p_2}$. We denote $q^*=R_{p_1}+R_{p_2}$. Then, using the fact that the $q^*$ first moments of
$\Lambda_{j_1,j_2}^{p_1,p_2}$ are null,
\begin{equation}
\textrm{Cov}\left[W_{j_1}^{p_1}(k_1),W_{j_2}^{p_2}(k_2)\right] \sim ~ C_1|\alpha|^{2d-1-q^*}+R_{q^*+1},
\end{equation}
with $C_1$ a bounded constant, and:
\begin{eqnarray}
  \abs{R_{q^*+1}} &=& \abs{\cos(\nu\alpha)}|\alpha|^{2d-1}\abs{\sum_{i=q^*+1}^{\infty}{2d-1\choose
i}\int\left(\frac{h}{\alpha}\right)^{i}\Lambda_{j_1,j_2}^{p_1,p_2}(h)dh} \nonumber\\
  &\leq& |\alpha|^{2d-1}{2d-1\choose q^*}\abs{\sum_{i=q^*+1}^{\infty}\int\int\left(\frac{t-h}{\alpha}\right)^{i}\psi_{j_1}^{p_1}(t)\psi_{j_2}^{p_2}(h)dtdh} \nonumber\\
  &=& |\alpha|^{2d-1}{2d-1\choose q^*}\int\int\left|\psi_{j_1}^{p_1}(t)\psi_{j_2}^{p_2}(h)\right|dtdh\sum_{i=q^*+1}^{\infty}\beta^{i} \nonumber\\
  &=& C_2|\alpha|^{2d-1}\sum_{i=1}^{\infty}\beta^{q^*+i}
  \leq C_3|\alpha|^{2d-1-q^*-1}
\end{eqnarray}
where $\beta=\sup_{t,h}\left|\frac{t-h}{\alpha}\right|$, $C_2$ and $C_3$ are finite constants. Finally,
\begin{equation}
\textrm{Cov}\left[W_{j_1}^{p_1}(k_1),W_{j_2}^{p_2}(k_2)\right]=O\left(|2^{j_1}k_1-2^{j_2}k_2|^{2d-1-q^*}\right),
\end{equation}
for $|2^{j_1}k_1-2^{j_2}k_2|>(N^*+1)(2^{j_1}+2^{j_2})$ and $q^*=R_{p_1}+R_{p_2}$.


\item If $p_1 \neq 0$ {\bf or} $p_2 \neq 0$ then $\mathcal{M}^{j_1,j_2}_{p_1,p_2}(m) = 0,$ for $0\leq
m<R_{\max(p_1,p_2)}$. Following the same steps as above, we prove the second statement of Theorem \ref{Cov_Wave} with $q^*=R_{\max(p_1,p_2)}$.


\item $p_1=p_2=0$: $\mathcal{M}^{j_1,j_2}_{p_1,p_2}(m) = 0$ for $1\leq m<2q.$ In this particular case, we have
necessarily $j_1=j_2=j$. We must then upper-bound the covariance. From (\ref{eq_general}), we have:
\begin{equation}
\textrm{Cov}\left[W_{j}^{0}(k_1),W_{j}^{0}(k_2)\right] \sim ~ C_0|\alpha|^{2d-1}+C_1|\alpha|^{2d-1-2q}+R_{2q+1},
\end{equation}
where
\[
C_0=\int\cos(\nu(h+\alpha))\Lambda_{j,j}^{0,0}(h)dh,\ \ \ \
C_1=\cos(\nu\alpha)\frac{(2d-1)!}{(2q)!(2d-1-2q)!}\int h^{2q}\Lambda_{j,j}^{0,0}(h)dh,
\]
and
\begin{eqnarray}
|R_{2q+1}| &\leq& |\cos(\nu\alpha)||\alpha|^{2d-1}|\sum_{i=2q+1}^{\infty}{2d-1\choose
i}\int\left(\frac{h}{\alpha}\right)^{i}\Lambda_{j,j}^{0,0}(h)dh| = O\left(|\alpha|^{2d-1-2q-1}\right).
\end{eqnarray}
As previously, when $\alpha$ is large:
\begin{eqnarray}
C_0 &\sim& \cos(\nu\alpha)\int\int\phi_j(t)\phi_j(t-h)dtdh = 2^j\cos(\nu\alpha)\left|\Phi(0)\right|^2=2^j\cos(\nu\alpha).
\end{eqnarray}
Finally, using a similar argument as in the previous cases, we find that for $|k_1-k_2|>2(N^*+1)$:
\begin{eqnarray}
\textrm{Cov}\left[W_{j}^{0}(k_1),W_{j}^{0}(k_2)\right] &=& O(|2^j(k_1-k_2)|^{2d-1}).
\end{eqnarray}
\end{enumerate}
\endproof

\bibliographystyle{elsart-num}
\bibliography{synthesisGG}

\begin{thebibliography}{10}
\expandafter\ifx\csname url\endcsname\relax
  \def\url#1{\texttt{#1}}\fi
\expandafter\ifx\csname urlprefix\endcsname\relax\def\urlprefix{URL }\fi

\bibitem{CoiWic92}
R.~Coifman, M.~Whickerhauser, Entropy-based algorithms for best basis
  selection, IEEE Transactions on Information Theory 38 (1992) 713--718.

\bibitem{Whi01}
B.~Whitcher, Simulating gaussian stationary processes with unbounded spectra,
  Journal of Computational and Graphical Statistics 10 (2001) 112--134.

\bibitem{Ber94}
J.~Beran, Statistics for Long-Memory Processes, Chapman and Hall, London, 1994.

\bibitem{Hosk84}
J.~Hosking, Modeling persistence in hydrological time series using fractional
  differencing, Water Resources Research 20 (1984) 1898--1908.

\bibitem{DavHar87}
R.~Davies, D.~Harte, Tests for hurst effect, Biometrika 74 (1987) 95--102.

\bibitem{Wor96}
G.~Wornell, Signal Processing with Fractals: A Wavelet-Based Approach, Prentice
  Hall, 1996.

\bibitem{MccWal96}
M.~McCoy, A.~Walden, Wavelet analysis and synthesis of stationary long-memory
  processes, Journal of Computational and Graphical Statistics 5 (1996) 1--31.

\bibitem{DerTew93}
M.~Deriche, A.~Tewfik, Maximum likelihood estimation of the parameters of
  discrete fractionally differenced gaussian noise process., IEEE Transactions
  on Signal Processing 41 (1993) 2977--2989.

\bibitem{TewKim92}
A.~Tewfik, M.~Kim, Correlation structure of the discrete wavelet coefficient of
  fractional brownian motion, IEEE Transactions on Information Theory 38 (1992)
  904--909.

\bibitem{Jen99}
M.~Jensen, Using wavelets to obtain a consistent ordinary least squares
  estimator of the long-memory parameter, Journal of Forecasting 18 (1999)
  17--32.

\bibitem{PerWal00}
D.~Percival, A.~Walden, Wavelet Methods for Time Series Analysis, Cambridge:
  Cambridge University Press, 2000.

\bibitem{Wor90}
G.~W. Wornell, A {Karhunen-Lo\`eve}-like expansion for 1/f process via
  wavelets, IEEE Trans. Inf. Theory 36~(4) (1990) 859--861.

\bibitem{Flan92}
P.~Flandrin, Wavelet analysis and synthesis of fractional brownian motion, IEEE
  Trans. Inf. Theory 38 (1992) 910--917.

\bibitem{Dijk94}
P.~Dijkerman, R.~Mazumdar, On the correlation structure of the wavelet
  coefficients of fractional brownian motion, IEEE Trans. Inf. Theory 45 (1994)
  1609--1612.

\bibitem{Jen00}
M.~Jensen, An alternative maximum likelihood estimator of long-memory processes
  using compactly supported wavelets, Journal of Economic Dynamics and Control
  24 (2000) 361--387.

\bibitem{MalZhaPap98}
S.~Mallat, Z.~Zhang, G.~Papanicolaou, Adaptive covariance estimation of locally
  stationary processes, Annals of Statistics 26 (1998) 1--47.

\bibitem{DonMalVon98}
D.~Donoho, S.~Mallat, R.~von Sachs, Estimating covariances of locally
  stationary processes: rates of convergence of best basis methods, Tech. rep.,
  University of Berkeley (1998).

\bibitem{Donoho97}
D.~Donoho, Cart and best-ortho-basis: a connection, Ann. Statist. 25~(5) (1989)
  1870--1911.

\bibitem{CoiMeyWic92}
R.~Coifman, Y.~Meyer, M.~Whickerhauser, Wavelet analysis and signal
  processing,, In Wavelets and their applications, Boston. (1992) 153--178.

\bibitem{Wic94}
M.~Wickerhauser, Adapted Wavelet Analysis from Theory to Software, IEEE Press,
  1994.

\bibitem{Mal98}
S.~Mallat, A wavelet tour in signal processing, Academic Press, 1998.

\bibitem{GraZhaWoo89}
H.~Gray, N.~Zhang, W.~Woodward, On generalized fractional processes, JTSA 10
  (1989) 233--257.

\bibitem{Gray98}
W.~Woodward, Q.~Cheng, H.~Gray, Long memory time series, J. of Time Series
  Analysis 9~(4) (1998) 485--504.

\bibitem{Hos81}
J.~Hosking, Fractional differencing, Biometrika 68 (1981) 165--176.

\bibitem{GraJoy80}
C.~Granger, R.~Joyeux, An introduction to long-memory time series models and
  fractional differencing, JTSA 1 (1980) 15--29.

\bibitem{Chu96}
C.~Chung, Estimating a generalized long-memory process, Journal of Econometrics
  73 (1996) 237--259.

\bibitem{Wic91}
M.~Wickerhauser, Fast approximate factor analysis, in: Curves and Surfaces in
  Computer Vision and Grpahics II, Vol. 1610, SPIE, 1991, pp. 23--32.

\bibitem{And86}
J.~Andel, Long memory time series, Kybernetika 22 (1986) 105--123.

\bibitem{Chung96}
C.~Chung, A generalized integrated autoregressive moving average process, J. of
  Time Series Analysis 17 (1996) 111--140.

\bibitem{Lapsa97}
P.~Lapsa, Determination of gegenbauer-type random process models, Signal
  Processing 63 (1997) 73--90.

\bibitem{Craigmile04}
P.~Craigmile, P.~Guttorp, D.~Percival, Wavelet-based parameter estimation for
  polynomial contaminated fractionally differenced processes, IEEE Trans. on
  Sig. Proc. 53~(8) (2004) 3151--3161.

\bibitem{Bey92}
G.~Beylkin, On the representation of operators in bases of compactly supported
  wavelets, SIAM Journal on Numerical Analysis 29~(6) (1992) 1716--1740.

\end{thebibliography}



\newpage
\vspace*{-1cm}

\begin{table}[H]
\begin{center}
\scriptsize
\hspace*{-1.5cm}
\begin{tabular}{p{4.5cm}p{4.1cm}p{4.1cm}p{4.1cm}}
$d=0.4,\ \nu=1/12,\ \lambda=20.7084$ & $d=0.2,\ \nu=1/12,\ \lambda=0.7428$ & $d=0.3,\ \nu=0.016,\ \lambda=10.0526$ & $d_1=d_2=0.3,\ \nu_1=1/40, \nu_2=1/5,\ \lambda= 6.0472$ \\
\begin{tabular}{p{0.2cm}p{1.5cm}p{1.25cm}} \hline
$q$   &  Whitcher basis  &  Our basis \\ \hline \hline \multicolumn{3}{p{4cm}}{\centering
Daubechies}
\\ \hline \hline
2  & 2728.6 & 1494.5  \\
4  & 1116.7 & 686.2   \\
6  & 750.4  & 441.8   \\
8  & 632.7  & 352.4   \\
10 & 421.7  & 308.2   \\ \hline \multicolumn{3}{p{4cm}}{\centering Symmlet} \\ \hline \hline
4  & 2211.0 & 677.0   \\
6  & 1371.0 & 444.7   \\
8  & 980.4  & 341.7   \\
10 & 693.1  & 297.3   \\ \hline \multicolumn{3}{p{4cm}}{\centering Coiflet} \\ \hline \hline
2  & 2697.3 & 1081.1  \\
4  & 1449.4 & 638.4   \\
6  & 841.1  & 412.7   \\
8  & 703.2  & 327.8   \\
10 & 581.2  & 287.7   \\ \hline \multicolumn{3}{p{4cm}}{\centering Battle-Lemari\'e} \\ \hline \hline
2  &   -    & 657.3   \\
4  &   -    & 267.5   \\
6  &   -    & 247.8   \\ \hline
\end{tabular} &
\begin{tabular}{p{1.5cm}p{1.5cm}} \hline
Whitcher basis  &  Our basis \\ \hline \hline
\multicolumn{2}{p{4cm}}{\centering Daubechies} \\ \hline \hline
105.1  & 52.3  \\
47.3   & 31.1  \\
31.9   & 23.3  \\
28.9   & 20.2  \\
21.9   & 18.4  \\ \hline \multicolumn{2}{p{4cm}}{\centering Symmlet} \\ \hline \hline
86.4   & 30.8  \\
54.2   & 23.4  \\
39.1   & 20.1  \\
28.2   & 18.2  \\ \hline \multicolumn{2}{p{4cm}}{\centering Coiflet} \\ \hline \hline
106.0  & 42.3  \\
58.9   & 29.6  \\
35.0   & 22.4  \\
29.0   & 19.5  \\
26.6   & 17.6  \\ \hline \multicolumn{2}{p{4cm}}{\centering Battle-Lemari\'e} \\ \hline \hline
  -    & 29.4  \\
  -    & 16.1  \\
  -    & 14.4  \\ \hline
\end{tabular} &  
\begin{tabular}{p{1.5cm}p{1.5cm}} \hline
Whitcher basis  &  Our basis \\ \hline \hline
\multicolumn{2}{p{4cm}}{\centering Daubechies} \\ \hline \hline
230.6  & 141.4 \\
188.0  & 124.7 \\
141.8  & 116.7 \\
144.0  & 118.8 \\
140.3  & 115.0 \\ \hline \multicolumn{2}{p{4cm}}{\centering Symmlet} \\ \hline \hline
214.4  & 120.9 \\
210.4  & 116.2 \\
178.6  & 114.7 \\
138.8  & 113.5 \\ \hline \multicolumn{2}{p{4cm}}{\centering Coiflet} \\ \hline \hline
224.8  & 132.7 \\
215.2  & 121.0 \\
141.2  & 116.1 \\
140.1  & 114.9 \\
138.6  & 113.2 \\\hline \multicolumn{2}{p{4cm}}{\centering Battle-Lemari\'e} \\ \hline \hline
  -    & 121.3 \\
  -    & 111.8 \\
  -    & 110.7 \\ \hline
\end{tabular} &
\begin{tabular}{p{1.5cm}p{1.5cm}} \hline
Whitcher basis  &  Our basis \\ \hline \hline
\multicolumn{2}{p{4cm}}{\centering Daubechies} \\ \hline \hline
958.5  & 639.1 \\
511.4  & 515.7 \\
444.1  & 422.5 \\
408.0  & 361.9 \\
400.5  & 315.3 \\ \hline \multicolumn{2}{p{4cm}}{\centering Symmlet} \\ \hline \hline
505.3  & 513.8 \\
445.0  & 423.5 \\
409.0  & 359.7 \\
401.2  & 313.7 \\ \hline \multicolumn{2}{p{4cm}}{\centering Coiflet} \\ \hline \hline
713.5  & 589.5 \\
467.1  & 500.7 \\
416.1  & 406.2 \\
363.3  & 342.7 \\
349.9  & 297.8 \\ \hline \multicolumn{2}{p{4cm}}{\centering Battle-Lemari\'e} \\ \hline \hline
  -    & 374.6 \\
  -    & 238.1 \\
  -    & 191.1 \\ \hline
\end{tabular}
\end{tabular}
\end{center}
\caption{$S$ score as a function of the number of vanishing moments for each wavelet family.}
\label{table:exact}
\end{table}

\vspace*{-1cm}

\begin{table}
\vspace*{-1.5cm}
\scriptsize
\begin{center}
\begin{tabular}{cc}
$d=0.4,\ \nu=1/12,\ \lambda=20.7084$ & $d=0.2,\ \nu=1/12,\ \lambda=0.7428$ \\
\begin{tabular}{@{}p{0.2cm}p{1.25cm}@{}@{}p{1.25cm}@{}@{}p{1.25cm}@{}@{}p{1.25cm}} \hline
\multicolumn{5}{p{5.2cm}}{\centering $B_{\text{Hosk}}=277.6$} \\ \hline
$q$   &  $B^W$ & $B_{\text{pen}}^W$ & $B^{CF}$ & $B_{\text{pen}}^{CF}$ \\ \hline \hline \multicolumn{5}{p{5.2cm}}{\centering Daubechies} \\ \hline \hline
2  & 2753.3   &   6977.8   &   2785.0	&   2992.1\\
4  & 1517.6   &   3629.8   &   1578.2	&   1785.2\\
6  & 1410.6   &   2777.4   &   1177.6	&   1384.7\\
8  &  929.7   &   1840.8   &    996.4	&   1203.5\\
10 & 1052.9   &   1694.9   &    784.6	&    991.6\\ \hline \multicolumn{5}{p{5.2cm}}{\centering Symmlet} \\ \hline \hline
4  & 1700.4   &   3812.6   &	1590.9   &   1797.9\\
6  & 1035.6   &   2402.4   &	1225.3   &   1432.4\\
8  & 1180.8   &   2091.9   &	1079.6   &   1286.7\\
10 &  718.9   &   1360.9   &	 797.4	 &   1004.5\\ \hline \multicolumn{5}{p{5.2cm}}{\centering Coiflet} \\ \hline \hline
2  & 2194.3   &   5238.4   &	2336.7   &   2543.8\\
4  & 1751.5   &   3118.3   &	1696.7   &   1903.8\\
6  & 1213.8   &   1855.8   &	1154.0   &   1361.1\\
8  & 1138.9   &   1656.7   &	1038.4   &   1245.5\\
10 &  995.2   &   1409.3   &	 986.3   &   1193.4\\ \hline \multicolumn{5}{p{5.2cm}}{\centering Battle-Lemari\'e} \\ \hline \hline
2  &   -     &	   -      &    1588.8	&   1795.9\\
4  &   -     &	   -      &     995.4	&   1202.4\\
6  &   -     &	   -      &     867.0	&   1074.0\\ \hline
\end{tabular} &
\begin{tabular}{p{1.25cm}@{}@{}p{1.25cm}@{}@{}p{1.25cm}@{}@{}p{1.25cm}} \hline
\multicolumn{4}{p{5cm}}{\centering $B_{\text{Hosk}}=1.72$} \\ \hline
$B^W$ & $B_{\text{pen}}^W$ & $B^{CF}$ & $B_{\text{pen}}^{CF}$ \\ \hline \hline \multicolumn{4}{p{5cm}}{\centering Daubechies} \\ \hline \hline
41.4   &    192.9   &	36.6   &     44.0\\
15.6   &     91.3   &	13.0   &     20.4\\
14.3   &     63.3   &	10.8   &     18.2\\
10.9   &     43.5   &	 6.9   &     14.3\\
12.6   &     35.6   &	 4.7   &     12.2\\\hline \multicolumn{4}{p{5cm}}{\centering Symmlet} \\ \hline \hline
25.5   &    101.3   &	 16.2	&     23.7\\
10.0   &     59.1   &	 11.6	&     19.0\\
 8.8   &     41.5   &	 11.6	&     19.0\\
 8.7   &     31.7   &	 12.6	&     20.0\\ \hline \multicolumn{4}{p{5cm}}{\centering Coiflet} \\ \hline \hline
27.8   &    137.0   &  34.8   &     42.2\\
18.3   &     67.3   &  18.9   &     26.3\\
10.2   &     33.2   &  10.8   &     18.3\\
 5.1   &     23.7   &   9.2   &     16.6\\
11.8   &     26.7   &   6.4   &     13.8\\ \hline \multicolumn{4}{p{5cm}}{\centering Battle-Lemari\'e} \\ \hline \hline
-   &	 -   &  18.0   &     25.4\\
-   &	 -   &   9.0   &     16.4\\
-   &	 -   &   5.6   &     13.0\\ \hline
\end{tabular} \\\\
$d=0.3,\ \nu=0.016,\ \lambda=10.0526$ & 
$d_1=d_2=0.3,\ \nu_1=1/40, \nu_2=1/5,\ \lambda= 6.0472$ \\
\begin{tabular}{p{1.25cm}@{}@{}p{1.25cm}@{}@{}p{1.25cm}@{}@{}p{1.25cm}} \hline
\multicolumn{4}{p{5cm}}{\centering $B_{\text{Hosk}}=34.7$} \\ \hline
$B^W$ & $B_{\text{pen}}^W$ & $B^{CF}$ & $B_{\text{pen}}^{CF}$ \\ \hline \hline \multicolumn{4}{p{5cm}}{\centering Daubechies} \\ \hline \hline
499.6	&    751.0   &  280.0	&    380.6\\
258.9	&    490.1   &  231.7	&    332.3\\
299.5	&    440.3   &  246.9	&    347.5\\
243.9	&    384.7   &  185.4	&    285.9\\
310.6	&    451.3   &  345.3	&    445.8\\\hline \multicolumn{4}{p{5cm}}{\centering Symmlet} \\ \hline \hline
287.1	&    518.3   &  284.0	&    384.5\\
232.6	&    373.3   &  312.2	&    412.7\\
362.3	&    503.0   &  159.6	&    260.2\\
218.4	&    359.2   &  225.2	&    325.7\\  \hline \multicolumn{4}{p{5cm}}{\centering Coiflet} \\ \hline \hline
280.4	&    511.7   &  293.6	&    394.1\\
208.0	&    348.7   &  166.0	&    266.5\\
377.5	&    518.3   &  268.5	&    369.0\\
273.5	&    414.2   &  198.6	&    299.1\\
257.5	&    398.3   &  292.8	&    393.4\\ \hline \multicolumn{4}{p{5cm}}{\centering Battle-Lemari\'e} \\ \hline \hline
 -   &    -   &  308.7   &    409.2\\
 -   &    -   &  321.5   &    422.0\\
 -   &    -   &  324.4   &    424.9\\ \hline
\end{tabular} &
\begin{tabular}{p{1.25cm}@{}@{}p{1.25cm}@{}@{}p{1.25cm}@{}@{}p{1.25cm}} \hline
\multicolumn{4}{p{5cm}}{\centering $B_{\text{Hosk}}=44.3$} \\ \hline
$B^W$ & $B_{\text{pen}}^W$ & $B^{CF}$ & $B_{\text{pen}}^{CF}$ \\ \hline \multicolumn{4}{p{5cm}}{\centering Daubechies} \\ \hline \hline
382.1	&   1639.9   &   386.3   &    489.1\\
301.9	&    906.6   &   291.7   &    394.5\\
199.6	&    707.6   &   161.1   &    263.9\\
157.1	&    628.7   &   217.2   &    320.0\\
205.4	&    677.1   &   215.5   &    318.3\\\hline \multicolumn{4}{p{5cm}}{\centering Symmlet} \\ \hline \hline
274.9	&    879.6   &  282.7	&    385.5\\
247.1	&    755.0   &  264.7	&    367.5\\
202.5	&    674.2   &  174.7	&    277.5\\
172.1	&    643.7   &  155.5	&    258.3\\ \hline \multicolumn{4}{p{5cm}}{\centering Coiflet} \\ \hline \hline
381.7	&   1288.8   &  382.0	&    484.8\\
284.8	&    792.7   &  305.7	&    408.5\\
225.9	&    697.5   &  217.8	&    320.6\\
127.6	&    526.7   &  218.6	&    321.4\\
121.5	&    514.6   &  207.6	&    310.4\\ \hline \multicolumn{4}{p{5cm}}{\centering Battle-Lemari\'e} \\ \hline \hline
 -   &    -   &  316.6   &    419.4\\
 -   &    -   &  175.8   &    278.6\\
 -   &    -   &  142.0   &    244.8\\ \hline
\end{tabular}
\end{tabular}
\end{center}
\caption{Squared difference ($B$) and its penalized version ($B_{\text{pen}}$) for our basis (superscript CF) and Whitcher basis \cite{Whi01} (superscript W) as a function of the number of vanishing moments for each wavelet family.}
\label{table:simu} 
\end{table}

\end{document}